\documentclass[11pt]{article}
\usepackage{amscd}
\usepackage{amssymb,amsfonts,amsmath,psfrag,amscd,cite}
\usepackage{amsmath}
  \usepackage{paralist}
  \usepackage{graphics}
  \usepackage{epsfig}
  \usepackage[colorlinks=true]{hyperref}
  \usepackage[all]{xy}
\newtheorem{theo}{Theorem}

\newtheorem{defi}[theo]{Definition}

\newtheorem{rema}[theo]{Remark}

\makeatletter \@addtoreset{equation}{section}
\@addtoreset{theo}{section} \makeatother

\begin{document}


\date{}
\title{Regular Reduction of Controlled Magnetic Hamiltonian System
with Symmetry of the Heisenberg Group}
\author{Hong Wang \\
School of Mathematical Sciences and LPMC,\\
Nankai University, Tianjin 300071, P.R.China\\
E-mail: hongwang@nankai.edu.cn}
\date{\emph{Dedicated to Professor Tudor S. Ratiu
on the occasion of his 65th birthday}\\
July 7, 2018 } \maketitle

{\bf Abstract:} A controlled magnetic Hamiltonian (CMH) system is a regular controlled
Hamiltonian (RCH) system with magnetic symplectic form,
it is an important special case of RCH system. Note that there
is a magnetic term on the cotangent bundle of the Heisenberg group,
such that we can define a CMH system with symmetry of the
Heisenberg group. Since the set of the CMH systems with symmetries
is a subset of the RCH systems with symmetries, and it is not complete
under the regular point reduction of RCH system, in this paper, then we give
the regular point reduction of a CMH system with symmetry of the
Heisenberg group, and discuss the M-CH-equivalence and MR-CH-equivalence,
and prove the regular point reduction theorem for such CMH system.
In particular, we deduce the regular point reduced CMH
system on the generalization of coadjoint orbit of the Heisenberg
group by calculation in detail. As an application, we consider the
motion of the Heisenberg particle in a magnetic field.\\




\tableofcontents

\section{Introduction}

The reduction theory is an important subject and it is widely
studied in the theory of mathematics and mechanics,
as well as applications. In particular, the reduction theory
for mechanical system with symmetry has its
origin in the classical work of Euler, Lagrange, Hamilton, Jacobi,
Routh, Liouville and Poincar\'{e}, and its modern geometric
formulation in the general context of symplectic manifolds and
equivariant momentum maps is developed by Meyer, Marsden and
Weinstein; see Abraham and Marsden \cite{abma78} or Marsden and
Weinstein \cite{mawe74} and Meyer \cite{me73}. The main goal of
reduction theory in mechanics is to use conservation laws and the
associated symmetries to reduce the number of dimensions of a
mechanical system required to be described. So, such reduction
theory is regarded as a useful tool for simplifying and studying
concrete mechanical systems, and great developments have been
obtained around the Marsden-Weinstein reduction
in the theoretical study and applications
of mathematics, mechanics and physics. Also see Abraham et al.
\cite{abmara88}, Arnold \cite{ar89}, Libermann and Marle
\cite{lima87}, Marsden \cite{ma92}, Marsden et al.
\cite{mamiorpera07}, Marsden and Ratiu
\cite{mara99} and Ortega and Ratiu \cite{orra04} for more details.\\

The reduction theory of a Hamiltonian system with symmetry
on the cotangent bundle of a configuration manifold
is a very important special case of general reduction theory. We
first give a precise analysis for the geometrical structure of
Marsden-Weinstein symplectic reduced space of a Hamiltonian system
with symmetry on a cotangent bundle.\\

Let $Q$ be a smooth
manifold and $TQ$ the tangent bundle, $T^* Q$ the cotangent bundle
with a canonical symplectic form $\omega$. Assume that $\Phi:
G\times Q \rightarrow Q$ is a left smooth action of a Lie group $G$
on the manifold $Q$. The cotangent lift is the action of $G$ on
$T^\ast Q$, $\Phi^{T^\ast}:G\times T^\ast Q\rightarrow T^\ast Q$
given by $g\cdot \alpha_q=(T\Phi_{g^{-1}})^\ast\cdot
\alpha_q,\;\forall\;\alpha_q\in T^\ast_qQ,\; q\in Q$. The cotangent
lift of any proper (resp. free) $G$-action is proper (resp. free).
Assume that the cotangent lift action is symplectic with respect to the
canonical symplectic form $\omega$, and has an
$\operatorname{Ad}^\ast$-equivariant momentum map $\mathbf{J}:T^\ast
Q\to \mathfrak{g}^\ast$ given by
$<\mathbf{J}(\alpha_q),\xi>=\alpha_q(\xi_Q(q)), $ where $\xi\in
\mathfrak{g}$, $\xi_Q(q)$ is the value of the infinitesimal
generator $\xi_Q$ of the $G$-action at $q\in Q$, $<,>:
\mathfrak{g}^\ast \times \mathfrak{g}\rightarrow \mathbb{R}$ is the
duality pairing between the dual $\mathfrak{g}^\ast $ and
$\mathfrak{g}$. If $\mu\in \mathfrak{g}^\ast $ is a regular
value of the momentum map $\mathbf{J}$, and $G_\mu=\{g\in
G|\operatorname{Ad}_g^\ast \mu=\mu \}$ is the isotropy subgroup of
the coadjoint $G$-action at the point $\mu$. From the
Marsden-Weinstein reduction, we know that the reduced space
$((T^*Q)_\mu=\mathbf{J}^{-1}(\mu)/G_\mu, \omega_\mu)$ is a symplectic manifold
with the reduced symplectic form $\omega_\mu$ given by
$\pi_\mu^\ast \omega_\mu=i_\mu^\ast \omega $,
see Abraham and Marsden \cite{abma78} or Marsden and
Weinstein \cite{mawe74}.\\

Moreover, from Marsden et al. \cite{mamiorpera07}
and Marsden and Perlmutter \cite{mape00},
we know that the classification of symplectic
reduced space of a Hamiltonian system
with symmetry on the cotangent bundle $T^* Q$ as follows.
(1) If $\mu=0$, the symplectic reduced space of cotangent bundle
$T^\ast Q$ at $\mu=0$ is given by
$((T^\ast Q)_\mu, \omega_\mu)= (T^\ast(Q/G), \omega_0)$,
where $\omega_0$ is the canonical symplectic form of cotangent
bundle $T^\ast (Q/G)$. Thus, the symplectic reduced space $((T^\ast
Q)_\mu, \omega_\mu)$ at $\mu=0$ is a symplectic vector bundle.
(2) If $\mu\neq0$, and $G$ is Abelian, then $G_\mu=G$, in this case the
Marsden-Weinstein symplectic reduced space $((T^*Q)_\mu, \omega_\mu)$ is
symplectically diffeomorphic to symplectic vector bundle $(T^\ast
(Q/G), \omega_0-B_\mu)$, where $B_\mu$ is a magnetic term. (3) If
$\mu\neq0$, and $G$ is not Abelian and $G_\mu\neq G$, in this case
the Marsden-Weinstein symplectic reduced space $((T^*Q)_\mu,
\omega_\mu)$ is symplectically diffeomorphic to a symplectic fiber
bundle over $T^\ast (Q/G_\mu)$ with fiber to be the coadjoint orbit
$\mathcal{O}_\mu$, see the cotangent bundle reduction
theorem---bundle version in Marsden et al.
\cite{mamiorpera07}, also see Marsden and Perlmutter
\cite{mape00}.\\

Thus, from the above discussion, we know that
the symplectic reduced space on a cotangent bundle may not be a
cotangent bundle. Therefore, the symplectic reduced system of a
Hamiltonian system with symmetry defined on the cotangent bundle
$T^*Q$ may not be a Hamiltonian system on a cotangent bundle, that
is, the set of Hamiltonian systems with symmetries on the cotangent
bundle is not complete under the Marsden-Weinstein reduction.\\

On the other hand, in mechanics, the phase space of a Hamiltonian
system is very often the cotangent bundle $T^* Q$ of a configuration
manifold $Q$. Therefore, it is a serious problem that the
Marsden-Weinstein reduction is not complete. For example, if we define
directly a controlled Hamiltonian system with symmetry on a
cotangent bundle, then it is possible that the Marsden-Weinstein
reduced system may not have definition.
Recently, in Marsden et al.\cite{mawazh10}, the authors set up the
regular reduction theory of regular controlled Hamiltonian
(RCH) systems on a symplectic fiber
bundle, by using momentum map and the associated reduced symplectic
forms and from the viewpoint of completeness of
Marsden-Weinstein symplectic reduction. In particular,
the authors introduce the notions of CH-equivalence, RpCH-equivalence and
RoCH-equivalence to emphasize explicitly the impact of external force
and control in the study of RCH systems. The research work
in Marsden et al.\cite{mawazh10} is very important, and
there are some generalizations around the work to have been done,
see Wang and Zhang \cite{wazh12}, Ratiu and Wang \cite{rawa12}
for more details. In addition, as the
applications, we also study the underwater vehicle-rotors system and
rigid spacecraft-rotors system, as well as Hamilton-Jacobi theory
for the RCH system with symmetry, which show the effect on regular
symplectic reduction of RCH system. See Wang \cite{wa17, wa13d,
wa13f, wa13e} for more details. These research works not only gave a
variety of reduction methods for controlled Hamiltonian systems, but
also showed a variety of relationships of the
controlled Hamiltonian equivalences of these systems.\\

The theory of controlled mechanical systems is a very important
subject, following the theoretical development of geometric
mechanics, a lot of important problems about this subject are being
explored and studied. We know that it is not easy to give the
precise analysis of geometrical and topological structures of the
phase spaces and the reduced phase spaces of various Hamiltonian
systems. The study of completeness of Hamiltonian reductions for
controlled Hamiltonian system with symmetry is related to the
geometrical structures of Lie group, configuration manifold and its
cotangent bundle, as well as the action ways of Lie group on the
configuration manifold and on its cotangent bundle. Our goal to do the
research is to set up the various perfect reduction theory for
controlled mechanical systems, along the ideas of Professor Jerrold
E. Marsden, by analyzing carefully the geometrical and topological
structures of the phase spaces of various mechanical systems,
see Wang \cite{wa18}.\\

Recently, we note that the Heisenberg group is an important Lie
group and it is a central extension of $\mathbb{R}^2$ by
$\mathbb{R}$, and hence it is also a motivating example for the
general theory of Hamiltonian reduction by stages, see Marsden et
al. \cite{mamiorpera07, mamipera98}, also see Capogna et al.
\cite{cadapaty07} and Montgomery \cite{mo02} for more details of the
geometry of the Heisenberg group. In particular, we note that there
is a magnetic term on the cotangent bundle of the Heisenberg group
$\mathcal{H}$, which is related to a curvature two-form of a
mechanical connection determined by the reduction of center action
of the Heisenberg group $\mathcal{H}, $ see Theorem 3.2 in $\S 3$,
such that we can define a kind of controlled magnetic Hamiltonian (CMH) system
with symmetry of the Heisenberg group $\mathcal{H}, $
and study the regular point reduction of such CMH system, by using the
reduction of the magnetic cotangent bundle of the Heisenberg group
$\mathcal{H} $ and from the viewpoint of completeness of
the regular point reduction of RCH system,
and discuss the magnetic reducible controlled
Hamiltonian (MR-CH) equivalence. These are the main works in this
paper. It is worthy of noting that the regular point reduction of a
CMH system with symmetry of the Heisenberg group $\mathcal{H} $
describes the impact of the special structure of symmetric group of the
CMH system for the regular reduction of a RCH system.
It is different from the regular point reduction of a RCH system
defined on a cotangent bundle with the canonical structure,
the regular point reduction of a CMH system can reveal
the deeper relationship of the intrinsic geometrical structure
of phase space of the RCH system.\\

A brief of outline of this paper is as follows. In the second
section, we review some relevant definitions and basic facts about
the Heisenberg group $\mathcal{H}, $ coadjoint $\mathcal{H}$-action
and coadjoint orbit, which will be used in subsequent sections. The
reduction of the magnetic cotangent bundle of the Heisenberg group
$\mathcal{H} $, and the magnetic term in the magnetic symplectic
form of cotangent bundle of the Heisenberg group $\mathcal{H}$,
which is related to a curvature two-form of a mechanical connection
determined by the reduction of center action of the Heisenberg group
$\mathcal{H}, $ are introduced in the third section. In the fourth
section, we introduce briefly some relevant definitions and basic
facts about the RCH systems defined on a symplectic fiber bundle and
on the cotangent bundle of a configuration manifold, respectively,
and RCH-equivalence, the regular point reducible RCH system with
symmetry, as well as CMH system, M-CH-equivalence, and
the regular point reducible CMH system with symmetry.
Even if a CMH system is also a RCH system, but the set of
CMH systems with symmetries is not a complete subset of the set of
RCH systems with symmetries under the regular point reduction of RCH system.
In the fifth section, we state that the CMH
system with symmetry of the Heisenberg group $\mathcal{H} $
is a regular point reducible CMH system, and
give its regular point reduced CMH
system on the generalization of coadjoint orbit of the
Heisenberg group $\mathcal{H} $ by calculation in detail.
Moreover, we discuss the MR-CH-equivalence for the
regular point reducible CMH system with symmetry of
the Heisenberg group $\mathcal{H} $, and prove the regular point reduction
theorem for such system, which explains the relationship between
MR-CH-equivalence for the regular point reducible CMH system with
symmetry of the Heisenberg group and M-CH-equivalence for the
associated regular point reduced CMH system. As an application of
the theoretical results, in the sixth section, we consider the
motion of the Heisenberg particle in a magnetic field, and we also
consider the magnetic term from the viewpoint of Kaluza-Klein
construction. These research works develop the theory of regular
reduction for the CMH system with symmetry of the Heisenberg group
and make us have much deeper understanding and recognition for the
geometrical structures of phase spaces
of controlled Hamiltonian systems.

\section{The Heisenberg Group and Its Lie Algebra}

In this paper, our goal is to define a kind of CMH system with
symmetry of the Heisenberg group $\mathcal{H} $ and to prove
the regular point reduction theorem for such CMH system,
by using the magnetic term and the magnetic symplectic form
on the cotangent bundle of the Heisenberg group $\mathcal{H}.$
In order to do these, in this section,
we first review some relevant definitions and basic
facts about the Heisenberg group $\mathcal{H}, $ coadjoint
$\mathcal{H}$-action and coadjoint orbit, which will be used in
subsequent sections. We shall follow the notations and conventions
introduced in Marsden et al. \cite{mamiorpera07,
mamipera98}.\\

We consider the commutative group $\mathbb{R}^2$ with its standard
symplectic form $\omega$, which is the usual area form on the
Euclidean plane, that is,
\begin{equation}
\omega(u,v)=u_1v_2-u_2v_1,
\end{equation}
where $u=(u_1,u_2),\; v=(v_1,v_2) \in \mathbb{R}^2. $ Define the set
$\mathcal{H}= \mathbb{R}^2\oplus \mathbb{R}$ with group
multiplication
\begin{equation}
(u,\alpha)(v,\beta)=(u+v, \alpha+\beta+\frac{1}{2}\omega(u,v))
\end{equation}
where $u,\; v \in \mathbb{R}^2 $ and $\alpha, \; \beta \in
\mathbb{R}. $ It is readily verified that this operation defines a
Lie group, and its identity element is $(0,0)$ and the inverse of
$(u,\alpha)$ is given by $(u,\alpha)^{-1}=(-u,-\alpha). $ This group
is called the {\bf Heisenberg group}, which is an important Lie
group and it is isomorphic to the upper triangular $3\times 3$
matrices with ones on the diagonal, and the isomorphism is given by
\begin{align}
 (u,\alpha) \mapsto \left[\begin{array}{ccc} 1 &  u_1  &
\alpha+\frac{1}{2}u_1u_2 \\ 0 &  1 &  u_2 \\ 0 & 0 & 1 \end{array}
\right].
\end{align}

It is worthy of noting that each element $(0,\alpha)$ in $\mathcal{H}$
commutes with every other element of $\mathcal{H}$, and by using the
nondegeneracy of the symplectic form $\omega$, we know that every
element of $\mathcal{H}$ that commutes with all other elements of
$\mathcal{H}$ is of the form $(0,\alpha). $ Hence, the subgroup
$A=\{(0,\alpha)\in \mathcal{H}| \alpha\in \mathbb{R}\}$, consisting
of pairs $(0,\alpha)$ in $\mathcal{H}$, is the {\bf center} of
$\mathcal{H}$ and $A\cong \mathbb{R}. $ Thus, the Heisenberg group
$\mathcal{H}$ is the {\bf central extension} of $\mathbb{R}^2$ by
$\mathbb{R}$ and $B= \omega$ is its group two-cocycle,
see Marsden et al.\cite{mamiorpera07}.\\

The Lie algebra of the Heisenberg group $\mathcal{H}$ is
$\mathfrak{\eta}\cong \mathbb{R}^2\oplus \mathbb{R}. $ We identify
$\mathfrak{\eta}$ with $\mathbb{R}^3$ via the Euclidean inner
product. In the following we can calculate the Lie algebra bracket
on $\mathfrak{\eta}$. We know that the left and right translation on
$\mathcal{H}$ induce the left and right action of $\mathcal{H}$ on
itself. The conjugation action $I: \mathcal{H} \to \mathcal{H}$ is
given by
\begin{equation}
I_{(u,\alpha)}((v,\beta))=(u,\alpha)(v,\beta)(u,\alpha)^{-1}=
(u,\alpha)(v,\beta)(-u,-\alpha)= (v, \beta+\omega(u,v)),
\end{equation}
for any $(u,\alpha),\; (v,\beta) \in \mathcal{H}$, which is the
inner automorphism on $\mathcal{H}$. By differentiating the above
formula with respect to $(v,\beta)$, we can see that the operator
$\operatorname{Ad}: \mathcal{H}\times \mathfrak{\eta} \rightarrow
\mathfrak{\eta}, $ which induces an adjoint action of $\mathcal{H}$
on $\mathfrak{\eta}$, is given by
\begin{equation}
 \operatorname{Ad}_{(u,\alpha)}(Y,b)=(Y, b+\omega(u,Y)),
\end{equation}
where $(Y,b)\in \mathfrak{\eta}. $ By differentiating once more the
above formula with respect to $(u,\alpha)$, we can get that the
operator $\operatorname{ad}: \mathfrak{\eta}\times \mathfrak{\eta}
\rightarrow \mathfrak{\eta}, $ that is, the Lie bracket $[,]:
\mathfrak{\eta}\times \mathfrak{\eta} \rightarrow \mathfrak{\eta},$
is given by
\begin{equation}
 \operatorname{ad}_{(X,a)}(Y,b)= [(X,a),(Y,b)]=(0, \omega(X,Y)),
\end{equation}
with a Lie algebra two-cocycle $C(X,Y)= \omega(X,Y), $ where $(X,a),
\; (Y,b)\in \mathfrak{\eta}\cong \mathbb{R}^2\oplus \mathbb{R}. $\\

The dual of Lie algebra $\mathfrak{\eta}$ of the Heisenberg group
$\mathcal{H}$ is $\mathfrak{\eta}^*\cong \mathbb{R}^2\oplus
\mathbb{R}. $ We also identify $\mathfrak{\eta}^*$ with
$\mathbb{R}^3$ via the Euclidean inner product. Note that the
adjoint representation of the Heisenberg group $\mathcal{H}$ is
defined by
$\operatorname{Ad}_{(u,\alpha)}(Y,b)=T_{(0,0)}I_{(u,\alpha)} (Y,b)=
T_{(u,\alpha)^{-1}}L_{(u,\alpha)} \cdot T_{(0,0)}
R_{(u,\alpha)^{-1}}(Y,b): \mathcal{H} \times\mathfrak{\eta}\to
\mathfrak{\eta}$, then the coadjoint representation of the
Heisenberg group $\mathcal{H}, $ that is, $\operatorname{Ad}^\ast:
\mathcal{H} \times\mathfrak{\eta}^\ast\to \mathfrak{\eta}^\ast, $ is
defined by the following equation
\begin{equation}\langle
\operatorname{Ad}_{(u,\alpha)^{-1}}^\ast(\mu,\nu),(Y,b)\rangle
=\langle (\mu,\nu), \operatorname{Ad}_{(u,\alpha)}(Y,b)\rangle,
\end{equation} where $(u,\alpha) \in \mathcal{H}, \;
(\mu,\nu)\in\mathfrak{\eta}^\ast$, and $(Y,b) \in\mathfrak{\eta}$,
and $\langle,\rangle$ denotes the natural pairing between
$\mathfrak{\eta}^\ast$ and $\mathfrak{\eta}$. Moreover, the
coadjoint representation $\operatorname{Ad}^\ast: \mathcal{H}
\times\mathfrak{\eta}^\ast\to \mathfrak{\eta}^\ast$ induces a left
coadjoint action of the Heisenberg group $\mathcal{H}$ on
$\mathfrak{\eta}^\ast$, which is given by
\begin{equation}
\operatorname{Ad}_{(u,\alpha)^{-1}}^\ast(\mu,\nu)=
(\mu+\nu\mathbb{J}u, \nu),
\end{equation} where $(u,\alpha) \in
\mathcal{H}, \; (\mu,\nu)\in\mathfrak{\eta}^\ast$, and
$\mathbb{J}u=\mathbb{J}(u_1,u_2)=(u_2,-u_1)$ is the matrix of the
standard symplectic form $\omega$ on $\mathbb{R}^2. $ Then the
coadjoint orbit $\mathcal{O}_{(\mu,\nu)}$ of this
$\mathcal{H}$-action through $(\mu,\nu)\in\mathfrak{\eta}^\ast$ are
that (1) $\mathcal{O}_{(\mu,0)}= \{(\mu,0)\}, $ and (2)
$\mathcal{O}_{(\mu,\nu\neq 0)}\cong \mathbb{R}^2 \times\{\nu\}. $
which are the immersed submanifolds of
$\mathfrak{\eta}^\ast$.\\

We know that $\mathfrak{\eta}^\ast$ is a Poisson manifold with
respect to the $(\pm)$ magnetic Lie-Poisson bracket
$\{\cdot,\cdot\}^B_\pm$ defined by
\begin{equation}
\{f,g\}^B_\pm(\mu,\nu):= \pm<(\mu,\nu), [\frac{\delta f}{
\delta(\mu,\nu)}, \frac{\delta g}{\delta (\mu,\nu)}]>-
\pi^*_{\mathcal{H}}B(0,0)(\frac{\delta f}{ \delta(\mu,\nu)},
\frac{\delta g}{\delta (\mu,\nu)}),
\end{equation}
for any $f,g\in C^\infty(\mathfrak{\eta}^\ast), $ and $(\mu,\nu)\in
\mathfrak{\eta}^\ast, $ where the element $\frac{\delta f}{\delta
(\mu,\nu)}\in\mathfrak{\eta}$ is defined by the equality
$$<(\rho,\tau),\frac{\delta f}{\delta (\mu,\nu)}>:=Df(\mu,\nu)\cdot
(\rho,\tau), $$ for any $(\rho,\tau)\in \mathfrak{\eta}^\ast$, see
Marsden and Ratiu \cite{mara99}. Thus, for the coadjoint orbit
$\mathcal{O}_{(\mu,\nu)}, \; (\mu,\nu)\in\mathfrak{\eta}^\ast$, the
magnetic orbit symplectic structures can be defined by
\begin{equation}\omega_{\mathcal{O}_{(\mu,\nu)}}^\pm(\rho,\tau)(\operatorname{ad}_{(X,a)}^\ast(\rho,\tau),
\operatorname{ad}_{(Y,b)}^\ast(\rho,\tau))=\pm
\langle(\rho,\tau),[(X,a),(Y,b)]\rangle-\pi^*_{\mathcal{H}}B(0,0)((X,a),(Y,b)),
\end{equation}
for any $(X,a),\;(Y,b)\in \mathfrak{\eta}, $ and $(\rho,\tau)\in
\mathcal{O}_{(\mu,\nu)}\subset \mathfrak{\eta}^\ast, $ which
coincide with the restriction of the magnetic Lie-Poisson brackets
on $\mathfrak{\eta}^\ast$ to the coadjoint orbit
$\mathcal{O}_{(\mu,\nu)}$. From the Symplectic Stratification
theorem we know that a finite dimensional Poisson manifold is the
disjoint union of its symplectic leaves, and its each symplectic
leaf is an injective immersed Poisson submanifold whose induced
Poisson structure is symplectic. In consequence, when
$\mathfrak{\eta}^\ast$ is endowed one of the magnetic Lie-Poisson
structures $\{\cdot,\cdot\}^B_\pm$, the symplectic leaves of the
Poisson manifolds $(\mathfrak{\eta}^\ast,\{\cdot,\cdot\}^B_\pm)$
coincide with the connected components of the magnetic orbits of the
elements in $\mathfrak{\eta}^\ast$ under the coadjoint action.

\section{Magnetic Cotangent Bundle Reduction}

In this section, we consider the reduction of magnetic cotangent
bundle of the Heisenberg group $T^* \mathcal{H}$ with magnetic
symplectic form $\omega_B= \omega_0-\pi^*_{\mathcal{H}}B,$ where
$\omega_0$ is the usual canonical symplectic form on $T^*
\mathcal{H}$, and $B$ is a closed two-form on $\mathcal{H}$, and
$\pi^*_{\mathcal{H}} B$ is the magnetic term on $T^*\mathcal{H}$,
the map $\pi_{\mathcal{H}}: T^* \mathcal{H}\rightarrow \mathcal{H}$
is the cotangent bundle projection and $\pi^*_{\mathcal{H}}: T^*
\mathcal{H} \rightarrow T^*T^*\mathcal{H}$. Defined the left
$\mathcal{H}$-action $\Phi: \mathcal{H}\times \mathcal{H}
\rightarrow \mathcal{H} $ given by
\begin{equation} \Phi((u,\alpha),(v,\beta)):=(u,\alpha)(v,\beta)= (u+v,
\alpha+\beta+\frac{1}{2}\omega(u,v)),
\end{equation} for any
$(u,\alpha),(v,\beta) \in \mathcal{H}, $ that is, the
$\mathcal{H}$-action on $\mathcal{H}$ is the left translation on
$\mathcal{H}$, which is free and proper, and leaves the two-form $B$
invariant. By using the local left trivialization of $T^* \mathcal{H}$, we
have that $T^\ast \mathcal{H}\cong \mathcal{H}\times
\mathfrak{\eta}^\ast$ (locally). We consider the cotangent lift of the
$\mathcal{H}$-action to the magnetic cotangent bundle $(T^\ast
\mathcal{H}, \omega_B)$, which is given by
\begin{equation}
\Phi^{T^*}: \mathcal{H}\times T^*\mathcal{H} \rightarrow
T^*\mathcal{H}, \;
\Phi^{T^*}((u,\alpha),((v,\beta),(\mu,\nu))):=((u,\alpha)(v,\beta),(\mu,\nu)),
\end{equation}
for any $(u,\alpha),(v,\beta) \in \mathcal{H},\; (\mu,\nu) \in
\mathfrak{\eta}^\ast, $ and it is also a free and proper action.\\

Assume that the $\mathcal{H}$-action on $(T^\ast \mathcal{H}, \omega_B)$ is
symplectic, and admits an $\operatorname{Ad}^\ast$-equivariant
momentum map $\mathbf{J}: T^\ast \mathcal{H}\rightarrow
\mathfrak{\eta}^\ast$. If
$(\mu,\nu)\in\mathfrak{\eta}^\ast$ is a regular value of
$\mathbf{J}$ and $G_{(\mu,\nu)}=\{(u,\alpha) \in
\mathcal{H}|\operatorname{Ad}_{(u,\alpha)}^\ast (\mu,\nu)=(\mu,\nu)
\}$ is the isotropy subgroup of coadjoint $\mathcal{H}$-action at
the point $(\mu,\nu)$. Since $G_{(\mu,\nu)} (\subset \mathcal{H})$
acts freely and properly on $\mathcal{H}$ and on $T^\ast
\mathcal{H}$, it follows that $G_{(\mu,\nu)}$ acts also freely and
properly on $\mathbf{J}^{-1}((\mu,\nu))$, so that the reduced space
$(T^\ast
\mathcal{H})_{(\mu,\nu)}=\mathbf{J}^{-1}((\mu,\nu))/G_{(\mu,\nu)}$
is a symplectic manifold with the reduced symplectic form
$\omega_{(\mu,\nu)}$ uniquely characterized by the relation
\begin{equation}\pi_{(\mu,\nu)}^\ast \omega_{(\mu,\nu)}=i_{(\mu,\nu)}^\ast
\omega_B. \label{3.3}\end{equation} The map
$i_{(\mu,\nu)}:\mathbf{J}^{-1}((\mu,\nu))\rightarrow T^\ast
\mathcal{H}$ is the inclusion and
$\pi_{(\mu,\nu)}:\mathbf{J}^{-1}((\mu,\nu))\rightarrow (T^\ast
\mathcal{H})_{(\mu,\nu)}$ is the projection. The pair $((T^\ast
\mathcal{H})_{(\mu,\nu)},\omega_{(\mu,\nu)})$ is called the regular
point reduced space of the magnetic cotangent bundle $(T^\ast
\mathcal{H},\omega_B)$ at
$(\mu,\nu)$.\\

If $(T^\ast \mathcal{H}, \omega_B)$ is a connected magnetic
symplectic manifold, and $\mathbf{J}: T^\ast \mathcal{H}\rightarrow
\mathfrak{\eta}^\ast$ is a non-equivariant momentum map with a
non-equivariance group one-cocycle $\sigma: \mathcal{H}\rightarrow
\mathfrak{\eta}^\ast$, which is defined by
$\sigma((u,\alpha)):=\mathbf{J}((u,\alpha)\cdot
z)-\operatorname{Ad}^\ast_{(u,\alpha)^{-1}}\mathbf{J}(z)$, where
$(u,\alpha)\in \mathcal{H}$ and $z\in T^\ast \mathcal{H}$. Then we
know that $\sigma$ produces a new affine action $\Theta:
\mathcal{H}\times \mathfrak{\eta}^\ast \rightarrow
\mathfrak{\eta}^\ast $ defined by
\begin{equation}
\Theta((u,\alpha),(\mu,\nu)):=\operatorname{Ad}^\ast_{(u,\alpha)^{-1}}(\mu,\nu)
+ \sigma((u,\alpha)),
\end{equation} where $(u,\alpha)\in
\mathcal{H}, \; (\mu,\nu) \in \mathfrak{\eta}^\ast$, with respect to
which the given momentum map $\mathbf{J}$ is equivariant. Since
$\mathcal{H}$ acts freely and properly on $T^\ast \mathcal{H}$, and
$\tilde{G}_{(\mu,\nu)}$ denotes the isotropy subgroup of $(\mu,\nu)
\in \mathfrak{\eta}^\ast$ relative to this affine action $\Theta$
and $(\mu,\nu)$ is a regular value of $\mathbf{J}$. Then the
quotient space $(T^\ast
\mathcal{H})_{(\mu,\nu)}=\mathbf{J}^{-1}((\mu,\nu))/\tilde{G}_{(\mu,\nu)}$
is also a symplectic manifold with the reduced symplectic form
$\omega_{(\mu,\nu)}$ uniquely characterized by (3.3), see
Ortega and Ratiu \cite{orra04}.\\

Moreover, from Abraham and Marsden \cite{abma78} and the above discussion
in $\S 2$, we can obtain the following theorem, which states that we
can describe the regular point reduced space of a magnetic cotangent
bundle by using the coadjoint orbit with the magnetic orbit
symplectic structure.
\begin{theo}
The coadjoint orbit $(\mathcal{O}_{(\mu,\nu)},
\omega_{\mathcal{O}_{(\mu,\nu)}}^{-}), \; (\mu,\nu)\in
\mathfrak{\eta}^\ast,$ is symplectically diffeomorphic to the
magnetic symplectic point reduced space $((T^\ast
\mathcal{H})_{(\mu,\nu)},\omega_{(\mu,\nu)})$ of the magnetic
cotangent bundle $(T^* \mathcal{H}, \omega_B)$.
\end{theo}

In the following we shall state that the magnetic term is related to
a curvature two-form of a mechanical connection, by the reduction of
center action of the Heisenberg group $\mathcal{H}. $ In $\S 2$, we
have known that the center of $\mathcal{H}$ is the subgroup
$A=\{(0,\alpha)\in \mathcal{H}| \alpha\in \mathbb{R}\}\cong
\mathbb{R}, $ and the Heisenberg group is the central extension of
$\mathbb{R}^2$ by $\mathbb{R}$, hence we can consider $\mathcal{H}$
as a right principal $\mathbb{R}$-bundle
$\mathcal{H}=\mathbb{R}^2\oplus \mathbb{R} \rightarrow \mathbb{R}^2,
$ and there is an induced right $\mathcal{H}$-invariant metric on
$\mathcal{H}$ as follows. In fact, we define that $<,>:
\mathfrak{\eta}\times \mathfrak{\eta} \rightarrow \mathbb{R}$ given
by $<(X,a),(Y,b)>= (X,Y)+ab, $ for any $(X,a), \; (Y,b)\in
\mathfrak{\eta}, $ where the Euclidean inner product $(,)$ in
$\mathbb{R}^2$ is used in the first summand and the multiplication
of real numbers in the second summand. For $(u,\alpha)\in
\mathcal{H}$ and $(X,a)\in T_{(u,\alpha)}\mathcal{H}, $ the tangent
of right translation on $\mathcal{H}$ is given by
\begin{equation}
T_{(u,\alpha)}R_{(v,\beta)}(X,a)=(X, a+ \frac{1}{2}\omega(X,v))\in
T_{(u,\alpha)(v,\beta)}\mathcal{H},
\end{equation}
and, in particular, we have that
\begin{equation}
T_{(u,\alpha)}R_{(u,\alpha)^{-1}}(X,a)=(X, a-
\frac{1}{2}\omega(X,u))\in \mathfrak{\eta}.
\end{equation}
Thus, the associated right $\mathcal{H}$-invariant metric on
$\mathcal{H}$ is given by
\begin{equation}
\ll(X,a),(Y,b)\gg_{(u,\alpha)}=
(X,Y)+ab-\frac{1}{2}a\omega(Y,u)-\frac{1}{2}b\omega(X,u)+\frac{1}{4}\omega(X,u)\omega(Y,u).
\end{equation}
Note that the exponential map of the Heisenberg group $exp:
\mathfrak{\eta} \rightarrow \mathcal{H}$ coincides with that of the
vector Lie group $(\mathbb{R}^3, +)$, that is, the identity map of
$\mathbb{R}^3. $ For a given $a \in \mathbb{R}, $ the infinitesimal
generator for the right $\mathbb{R}$-action of $\mathcal{H}$ on
$\mathcal{H}$ is given by
$$
a_{\mathcal{H}}(v,\beta)= \frac{d}{dt}(v,\beta)(0,ta)=(0,a).
$$
By combining these formulas and using the general formula for the
locked inertia tensor, see Marsden \cite{ma92}, we can get the
expression of the associated locked inertia tensor, that is,
\begin{equation}
<\mathbb{I}_{(u,\alpha)}(a),b>= \ll a_{\mathcal{H}}(u,\alpha),
b_{\mathcal{H}}(u,\alpha)\gg_{(u,\alpha)}=ab,
\end{equation} for any
$a,\; b\in \mathbb{R}. $ Moreover, for any $(X,a)\in
T_{(u,\alpha)}\mathcal{H}$ and $b\in \mathbb{R}, $ we have the
momentum map for the $\mathbb{R}$-action on $\mathcal{H}$ given by
\begin{equation}
<\mathbf{J}_{\mathbb{R}}(\ll(X,a),\cdot
\gg_{(u,\alpha)}),b>=\ll(X,a),(0,b)\gg_{(u,\alpha)}=(a-\frac{1}{2}\omega(X,u))b.
\end{equation} Thus, we can get the expression of the associated
mechanical connection $\mathcal{A}: T\mathcal{H}\rightarrow
\mathbb{R}$ given by
\begin{equation}
\mathcal{A}(u,\alpha)(X,a)=a-\frac{1}{2}\omega(X,u),
\end{equation}
and its exterior derivative is given by
\begin{equation}
\mathcal{B}=\mathbf{d}\mathcal{A}(u,\alpha)((X,a),(Y,b))=
\omega(X,Y),
\end{equation}
where $(X,a), \; (Y,b)\in T_{(u,\alpha)}\mathcal{H}, $ which offers
a Lie algebra valued, closed two-form $\mathcal{B}$ on
$\mathcal{H}$, that is, the curvature two-form of the mechanical
connection $\mathcal{A}$. For $\nu \in \mathbb{R}^* \cong
\mathbb{R}$, we can define the $\nu$-component of $\mathcal{B}$ by
$B^\nu=\nu \mathcal{B}, $ such that for any $(u,\alpha) \in
\mathcal{H}, \; (X,a), \; (Y,b) \in T_{(u,\alpha)}\mathcal{H}, \;
B^\nu(u,\alpha)((X,a),(Y,b))=\nu \mathcal{B}(u,\alpha)((X,a),(Y,b))=
\nu\omega(X,Y). $ This $B^\nu$ is an ordinary closed two-form on
$\mathcal{H}$, and $\pi^*_{\mathcal{H}} B^\nu$ is usually the
magnetic term on $T^* \mathcal{H}$. To sum up the above discussion,
we have the following theorem.
\begin{theo}
There is a magnetic term on the cotangent bundle of the Heisenberg
group $\mathcal{H}$, which is related to a curvature two-form of a
mechanical connection determined by the reduction of center action
of the Heisenberg group $\mathcal{H}. $
\end{theo}

\section{RCH System, CMH System and Regular Point Reduction}

In order to describe the regular point reduction of a CMH
system with symmetry of the Heisenberg group $\mathcal{H}$ and the
MR-CH-equivalence, in this section, we shall give some
relevant definitions and basic facts about RCH system,
RCH-equivalence, regular point reducible RCH
system, see Marsden et al.\cite{mawazh10},
as well as CMH system, M-CH-equivalence, regular point
reducible CMH system. We shall state that the set of
the CMH systems with symmetries is a subset of the set of
the RCH systems with symmetries,
and it is not complete under the regular point reduction
of RCH system. For convenience, we assume that all
controls appearing in this paper are the admissible controls.\\

At first, in order to describe
uniformly RCH systems defined on a cotangent bundle and on its
regular point reduced space, in terms of the completeness of
Marsden-Weinstein reduction, see Marsden et al.\cite{mawazh10},
we shall define a RCH system on a symplectic fiber bundle as follows.
Let $(E,M,N,\pi,G)$ be a fiber bundle and $(E, \omega_E)$ be a
symplectic fiber bundle. If for any function $H: E \rightarrow
\mathbb{R}$, we have a Hamiltonian vector field $X_H$ given by
the Hamilton's equation $\mathbf{i}_{X_H}\omega_E=\mathbf{d}H$,
then $(E, \omega_E, H )$ is a
Hamiltonian system. Moreover, if considering the external force and
control, we can define a kind of regular controlled Hamiltonian
(RCH) system on the symplectic fiber bundle $E$ as follows.

\begin{defi} (RCH System) A RCH system on $E$ is a 5-tuple
$(E, \omega_E, H, F, \mathcal{C})$, where $(E, \omega_E,\\ H )$ is a
Hamiltonian system, and the function $H: E \rightarrow \mathbb{R}$
is called the Hamiltonian, a fiber-preserving map $F: E\rightarrow
E$ is called the (external) force map, and a fiber submanifold
$\mathcal{C}$ of $E$ is called the control subset.
\end{defi}

Sometimes, $\mathcal{C}$ also denotes the set of fiber-preserving
maps from $E$ to $\mathcal{C}$. When a feedback control law $u:
E\rightarrow \mathcal{C}$ is chosen, the 5-tuple $(E, \omega_E, H,
F, u)$ denotes a closed-loop dynamic system. In particular, when $Q$
is a smooth manifold, and $T^\ast Q$ its cotangent bundle with a
symplectic form $\omega$ (not necessarily canonical symplectic
form), then $(T^\ast Q, \omega )$ is a symplectic vector bundle. If
we take that $E= T^* Q$, from above definition we can obtain a RCH
system on the cotangent bundle $T^\ast Q$, that is, 5-tuple $(T^\ast
Q, \omega, H, F, \mathcal{C})$.\\

In order to describe the dynamics of the RCH system
$(E,\omega_E,H,F,\mathcal{C})$ with a control law $u$, we can give a
good expression of the dynamical vector field of RCH system by using
the notations of vertical lifted maps of a vector along a fiber, see
Marsden et al.\cite{mawazh10}. In particular, in the case of
cotangent bundle, for a given RCH System $(T^\ast Q, \omega, H, F,
\mathcal{C})$, the dynamical vector field of the associated
Hamiltonian system $(T^\ast Q, \omega, H) $ is $X_H$, and it satisfies
Hamilton's equation $\mathbf{i}_{X_H}\omega=\mathbf{d}H$.
If considering the external force $F: T^*Q \rightarrow T^*Q, $
by using the notation of vertical lifted map of a vector along a fiber,
the change of $X_H$ under the action of $F$ is that
$$\mbox{vlift}(F)X_H(\alpha_x)
= \mbox{vlift}((TFX_H)(F(\alpha_x)), \alpha_x) =
(TFX_H)^v_\gamma(\alpha_x),$$ where $\alpha_x \in T^*_x Q, \; x\in Q
$ and $\gamma$ is a straight line in $T^*_x Q$ connecting
$F_x(\alpha_x)$ and $\alpha_x$. In the same way, when a feedback
control law $u: T^\ast Q \rightarrow \mathcal{C}$ is chosen, the
change of $X_H$ under the action of $u$ is that
$$\mbox{vlift}(u)X_H(\alpha_x)
= \mbox{vlift}((TuX_H)(u(\alpha_x)), \alpha_x) =
(TuX_H)^v_\gamma(\alpha_x).$$ In consequence, we can give an
expression of the dynamical vector field of RCH system as follows.
\begin{theo}
The dynamical vector field of a RCH system $(T^\ast
Q,\omega,H,F,\mathcal{C})$ with a control law $u$ is the synthetic
of Hamiltonian vector field $X_H$ and its changes under the actions
of the external force $F$ and control $u$, that is,
$$X_{(T^\ast Q,\omega,H,F,u)}(\alpha_x)= X_H(\alpha_x)+ \textnormal{vlift}(F)X_H(\alpha_x)
+ \textnormal{vlift}(u)X_H(\alpha_x),$$ for any $\alpha_x \in T^*_x
Q, \; x\in Q $. For convenience, it is simply written as
\begin{equation}X_{(T^\ast Q,\omega,H,F,u)}
=X_H +\textnormal{vlift}(F) +\textnormal{vlift}(u).
\end{equation}
\end{theo}
We also denote that $\mbox{vlift}(\mathcal{C})=
\bigcup\{\mbox{vlift}(u)X_H | \; u\in \mathcal{C}\}$. It is
worthy of noting that in order to deduce and calculate easily, we
always use the simple expression of dynamical vector
field $X_{(T^\ast Q,\omega,H,F,u)}$.\\

From the expression (4.1) of the dynamical vector
field of a RCH system, we know that under the actions of the external force $F$
and control $u$, in general, the dynamical vector
field is not Hamiltonian, and hence the RCH system is not
yet a Hamiltonian system. However,
it is a dynamical system closed relative to a
Hamiltonian system, and it can be explored and studied by extending
the methods for external force and control
in the study of Hamiltonian system. On the other hand,
we note that when a RCH system is given, the force map $F$ is
determined, but the feedback control law $u: T^\ast Q\rightarrow
\mathcal{C}$ could be chosen. In order to emphasize explicitly
the impact of external force and control in the study of the RCH systems,
by using the above expression of the dynamical vector field
of the RCH system, we can describe the feedback
control law how to modify the structure of a RCH system, and the controlled
Hamiltonian matching conditions and RCH-equivalence are induced as
follows.
\begin{defi}
(RCH-equivalence) Suppose that we have two RCH systems $(T^\ast
Q_i,\omega_i,H_i,F_i,\mathcal{C}_i),$ $ i= 1,2,$ we say them to be
RCH-equivalent, or simply, $(T^\ast
Q_1,\omega_1,H_1,F_1,\mathcal{C}_1)\stackrel{RCH}{\sim}\\ (T^\ast
Q_2,\omega_2,H_2,F_2,\mathcal{C}_2)$, if there exists a
diffeomorphism $\varphi: Q_1\rightarrow Q_2$, such that the
following controlled Hamiltonian matching conditions hold:

\noindent {\bf RCH-1:} The cotangent lifted map of $\varphi$, that
is, $\varphi^\ast= T^\ast \varphi:T^\ast Q_2\rightarrow T^\ast Q_1$
is symplectic, and $\mathcal{C}_1=\varphi^\ast (\mathcal{C}_2).$

\noindent {\bf RCH-2:} $Im[X_{H_1}+ \textnormal{vlift}(F_1)-
T\varphi^\ast X_{H_2}-\textnormal{vlift}(\varphi^\ast
F_2\varphi_\ast)]\subset \textnormal{vlift}(\mathcal{C}_1)$, where
the map $\varphi_\ast=(\varphi^{-1})^\ast: T^\ast Q_1\rightarrow
T^\ast Q_2$, and $T\varphi^\ast: TT^\ast Q_2\rightarrow TT^\ast
Q_1$, and $Im$ means the pointwise image of the map in brackets.
\end{defi}

In the following we consider the RCH system with symmetry and momentum map,
and give the regular point reducible RCH system.
Let $Q$ be a smooth manifold and $T^\ast Q$ its cotangent bundle
with the symplectic form $\omega$. Let $\Phi: G\times Q\rightarrow
Q$ be a smooth left action of a Lie group $G$ on $Q$, which is free
and proper. Then the cotangent lifted left action $\Phi^{T^\ast}:
G\times T^\ast Q\rightarrow T^\ast Q$ is also free and
proper. Assume that the action is symplectic and admits an
$\operatorname{Ad}^\ast$-equivariant momentum map $\mathbf{J}:T^\ast
Q\rightarrow \mathfrak{g}^\ast$, where $\mathfrak{g}$ is the Lie
algebra of $G$ and $\mathfrak{g}^\ast$ is the dual of
$\mathfrak{g}$. Let $\mu\in\mathfrak{g}^\ast$ be a regular value of
$\mathbf{J}$ and denote by $G_\mu$ the isotropy subgroup of the
coadjoint $G$-action at the point $\mu\in\mathfrak{g}^\ast$, which
is defined by $G_\mu=\{g\in G|\operatorname{Ad}_g^\ast \mu=\mu \}$.
Since $G_\mu (\subset G)$ acts freely and properly on $Q$ and on
$T^\ast Q$, then $Q_\mu=Q/G_\mu$ is a smooth manifold and that the
canonical projection $\rho_\mu:Q\rightarrow Q_\mu$ is a surjective
submersion. It follows that $G_\mu$ acts also freely and properly on
$\mathbf{J}^{-1}(\mu)$, so that the space $(T^\ast
Q)_\mu=\mathbf{J}^{-1}(\mu)/G_\mu$ is a symplectic manifold with the
symplectic form $\omega_\mu$ uniquely characterized by the relation
\begin{equation}\pi_\mu^\ast \omega_\mu=i_\mu^\ast
\omega. \label{6.1}\end{equation} The map
$i_\mu:\mathbf{J}^{-1}(\mu)\rightarrow T^\ast Q$ is the inclusion
and $\pi_\mu:\mathbf{J}^{-1}(\mu)\rightarrow (T^\ast Q)_\mu$ is the
projection. The pair $((T^\ast Q)_\mu,\omega_\mu)$ is called
Marsden-Weinstein reduced space of $(T^\ast Q,\omega)$ at $\mu$.\\

Assume that $H: T^\ast Q\rightarrow \mathbb{R}$ is a $G$-invariant
Hamiltonian, the flow $F_t$ of the Hamiltonian vector field $X_H$
leaves the connected components of $\mathbf{J}^{-1}(\mu)$ invariant
and commutes with the $G$-action, so it induces a flow $f_t^\mu$ on
$(T^\ast Q)_\mu$, defined by $f_t^\mu\cdot \pi_\mu=\pi_\mu \cdot
F_t\cdot i_\mu$, and the vector field $X_{h_\mu}$ generated by the
flow $f_t^\mu$ on $((T^\ast Q)_\mu,\omega_\mu)$ is Hamiltonian with
the associated Marsden-Weinstein reduced Hamiltonian function
$h_\mu:(T^\ast Q)_\mu\rightarrow \mathbb{R}$ defined by
$h_\mu\cdot\pi_\mu=H\cdot i_\mu$, and the Hamiltonian vector fields
$X_H$ and $X_{h_\mu}$ are $\pi_\mu$-related.\\

Since we can regard a Hamiltonian system on $T^*Q$
as a spacial case of a RCH system without external force and
control, then the set of Hamiltonian systems with symmetries
on $T^*Q$ is a subset of the set of RCH systems with symmetries on $T^*Q$.
If we first admit the Marsden-Weinstein reduction for a Hamiltonian system
with symmetry, then we may study the regular point
reduction for a RCH system with symmetry, as an extension of
Marsden-Weinstein reduction of a Hamiltonian system
under regular controlled Hamiltonian equivalence conditions.
On the other hand, note that the Marsden-Weinstein
reduced space $((T^*Q)_\mu, \omega_\mu)$ is symplectically
diffeomorphic to a symplectic fiber bundle. Thus, we can introduce a
regular point reducible RCH system as follows.
\begin{defi}
(Regular Point Reducible RCH System) A 6-tuple $(T^\ast Q, G,
\omega, H, F, \mathcal{C})$, where the Hamiltonian $H:T^\ast
Q\rightarrow \mathbb{R}$, the fiber-preserving map $F:T^\ast
Q\rightarrow T^\ast Q$ and the fiber submanifold $\mathcal{C}$ of\;
$T^\ast Q$ are all $G$-invariant, is called a regular point
reducible RCH system, if there exists a point
$\mu\in\mathfrak{g}^\ast$, which is a regular value of the momentum
map $\mathbf{J}$, such that the regular point reduced system, that
is, the 5-tuple $((T^\ast Q)_\mu,
\omega_\mu,h_\mu,f_\mu,\mathcal{C}_\mu)$, where $(T^\ast
Q)_\mu=\mathbf{J}^{-1}(\mu)/G_\mu$, $\pi_\mu^\ast
\omega_\mu=i_\mu^\ast\omega$, $h_\mu\cdot \pi_\mu=H\cdot i_\mu$,
$F(\mathbf{J}^{-1}(\mu))\subset \mathbf{J}^{-1}(\mu) $, $f_\mu\cdot
\pi_\mu=\pi_\mu \cdot F\cdot i_\mu$, $\mathcal{C}\cap
\mathbf{J}^{-1}(\mu)\neq \emptyset $, and
$\mathcal{C}_\mu=\pi_\mu(\mathcal{C}\cap \mathbf{J}^{-1}(\mu))$, is
a RCH system, which is simply written as $R_P$-reduced RCH system.
Here $((T^\ast Q)_\mu,\omega_\mu)$ is the $R_P$-reduced space, the function
$h_\mu:(T^\ast Q)_\mu\rightarrow \mathbb{R}$ is called the reduced
Hamiltonian, the fiber-preserving map $f_\mu:(T^\ast
Q)_\mu\rightarrow (T^\ast Q)_\mu$ is called the reduced (external)
force map, $\mathcal{C}_\mu$ is a fiber submanifold of \;$(T^\ast
Q)_\mu$ and is called the reduced control subset.
\end{defi}

In order to describe the impact of different geometric structures
of phase spaces of the RCH systems,
we introduce an important special case of RCH system,
that is, a controlled magnetic Hamiltonian (CMH) system.
We first give a magnetic Hamiltonian system on the
cotangent bundle $T^*Q$. Assume that
$T^*Q$ with the canonical symplectic form $\omega$,
and $B$ is a closed two-form on $Q$,
then $\omega_B= \omega- \pi_Q^*B$ is a symplectic form on $T^*Q$,
where $\pi_Q^*: T^*Q \rightarrow T^*T^*Q $. The
$\omega_B$ is called a magnetic symplectic form, and $\pi_Q^*B$ is called
a magnetic term on $T^*Q$, see Marsden et al. \cite{mamiorpera07}.
A magnetic Hamiltonian system is a 3-tuple $(T^\ast Q,\omega_B,H)$,
which is Hamiltonian system defined by a
magnetic symplectic form $\omega_B$. For a Hamiltonian $H$,
the dynamical vector field $X^B_H$, which is called
the magnetic Hamiltonian vector field,
satisfies the magnetic Hamilton's equation, that is,
$\mathbf{i}_{X^B_{H} }\omega_B= \mathbf{d}H $.\\

A controlled magnetic Hamiltonian (CMH) system on $T^*Q$ is
a 5-tuple $(T^\ast Q,\omega_B,H,F,\mathcal{C})$, which is a
magnetic Hamiltonian system $(T^\ast Q,\omega_B,H)$
with external force $F$ and control $\mathcal{C}$.
Thus, a CMH system is also a RCH system, but its symplectic structure
is given by a magnetic symplectic form, and the set of
the CMH systems is a subset of the set of the RCH systems. Moreover, we consider the set of
the CMH systems with symmetries, which is a subset of the set of
the RCH systems with symmetries, and we shall state that the subset
is not complete under the regular point reduction of RCH system.\\

In fact, from the classification of Marsden-Weinstein symplectic
reduced space of Hamiltonian system with symmetry
on the cotangent bundle $T^* Q$, we know that,
if the regular value of the momentum map $\mathbf{J}$,
$\mu\neq0$, and Lie group $G$ is Abelian, then the isotropy subgroup
$G_\mu=G$, in this case the
Marsden-Weinstein symplectic reduced space $((T^*Q)_\mu, \omega_\mu)$ is
symplectically diffeomorphic to symplectic vector bundle $(T^\ast
(Q/G), \omega_{B_\mu})$, where $\omega_{B_\mu}= \omega_0-B_\mu$ is a
magnetic symplectic form, and $B_\mu$ is a magnetic term, and
$\omega_0$ is the canonical symplectic form on $T^\ast (Q/G)$.
In fact, from Marsden and Ratiu \cite{mara99} we
know that this is a general phenomenon of momentum shifting
on the cotangent bundle $T^* Q$ under Marsden-Weinstein symplectic reduction.
In consequence, for a RCH system with symmetry
$(T^\ast Q, G, \omega, H, F, \mathcal{C})$, its regular point reduced
RCH system $(T^\ast (Q/G), \omega_{B_\mu}, h_\mu, f_\mu,\mathcal{C}_\mu)$
may be a CMH system. On the other hand, if a CMH system with symmetry
$(T^\ast Q, G, \omega_B, H, F, \mathcal{C})$ is regular point reducible,
when the regular value of the momentum map $\mathbf{J}$,
$\mu\neq0$, and Lie group $G$ is Abelian, in this case the
regular point reduced system is
$(T^\ast (Q/G), \omega_{\mu}, h_\mu, f_\mu,\mathcal{C}_\mu)$,
where the reduced symplectic form $\omega_{\mu}$ satisfies that
$\pi_\mu^\ast \omega_\mu=i_\mu^\ast \omega_B$, that is,
$\pi_\mu^\ast (\omega_0-B_\mu)=i_\mu^\ast (\omega- \pi_Q^*B)$, where
$\omega_{\mu}=\omega_0-B_\mu $, and $B_\mu$ is the magnetic term
for the reduced system.
In particular, when $\pi_\mu^\ast \cdot B_\mu=2i_\mu^\ast \cdot\pi_Q^*B$,
in this case, we have that
$\pi_\mu^\ast \cdot \omega_0=i_\mu^\ast (\omega+ \pi_Q^*B)$, that is,
from the viewpoint of the phenomenon of momentum shifting,
the regular point reduced system of the CMH system with symmetry
$(T^\ast Q, G, \omega+ \pi_Q^*B, H, F, \mathcal{C})$ is
$(T^\ast (Q/G), \omega_0, h_\mu, f_\mu,\mathcal{C}_\mu)$, which is
a RCH system with canonical symplectic form, that is,
a CMH system without magnetic term. Thus, the
reduced system of a CMH system with symmetry
may not be a CMH system under the regular point reduction
of RCH system, and hence the set of CMH systems with symmetries
is not complete under the regular point reduction
of RCH system. In order to give the regular point reduction
of CMH systems with symmetries and CH-equivalence for their reduced CMH systems,
we have to restrict the definitions of the RCH-equivalence
of RCH systems and of regular point reducible RCH system
to the set of CMH systems.\\

When we restrict the definition of RCH-equivalence of RCH systems to the set of
CMH systems, then we can describe the feedback control law to modify
the structures of CMH systems, and the M-CH-equivalence is induced as follows.
\begin{defi}
(M-CH-equivalence) Suppose that we have two CMH
systems $(T^\ast Q_i,\omega_{B_i},H_i,\\ F_i,\mathcal{C}_i),$ $ i=
1,2,$ we say them to be M-CH-equivalent, or simply, $(T^\ast
Q_1,\omega_{B_1},H_1,F_1,\mathcal{C}_1)\stackrel{M-CH}{\sim}\\
(T^\ast Q_2,\omega_{B_2},H_2,F_2,\mathcal{C}_2)$, if there exists a
diffeomorphism $\varphi: Q_1\rightarrow Q_2$, such that the
following magnetic controlled Hamiltonian matching conditions hold:

\noindent {\bf M-CH-1:} The cotangent lifted map of $\varphi$, that
is, $\varphi^\ast= T^\ast \varphi:T^\ast Q_2\rightarrow T^\ast Q_1$
is symplectic with respect to their magnetic symplectic forms, and
$\mathcal{C}_1=\varphi^\ast (\mathcal{C}_2).$

\noindent {\bf M-CH-2:} $Im[X^{B_1}_{H_1}+ \textnormal{vlift}(F_1)-
T\varphi^\ast X^{B_2}_{H_2}-\textnormal{vlift}(\varphi^\ast
F_2\varphi_\ast)]\subset \textnormal{vlift}(\mathcal{C}_1)$, where
the map $\varphi_\ast=(\varphi^{-1})^\ast: T^\ast Q_1\rightarrow
T^\ast Q_2$, and $T\varphi^\ast: TT^\ast Q_2\rightarrow TT^\ast
Q_1$, and $Im$ means the pointwise image of the map in brackets.
\end{defi}

Moreover, when we restrict the definition of regular point reducible RCH systems
with symmetries to the set of CMH systems with symmetries,
then we can give the definition of
regular point reducible CMH systems with symmetries as follows.
\begin{defi}
(Regular Point Reducible CMH System) A 6-tuple $(T^\ast Q,
G, \omega_B, H, F, \mathcal{C})$, where the Hamiltonian $H:T^\ast
Q\rightarrow \mathbb{R}$, the fiber-preserving map $F:T^\ast
Q\rightarrow T^\ast Q$ and the fiber submanifold $\mathcal{C}$ of\;
$T^\ast Q$ are all $G$-invariant, is called a regular point
reducible CMH system, if there exists a point
$\mu\in\mathfrak{g}^\ast$, which is a regular value of the momentum
map $\mathbf{J}$, such that the regular point reduced system, that
is, the 5-tuple $((T^\ast Q)_\mu,
\omega^B_\mu,h_\mu,f_\mu,\mathcal{C}_\mu)$, where $(T^\ast
Q)_\mu=\mathbf{J}^{-1}(\mu)/G_\mu$, $\pi_\mu^\ast
\omega^B_\mu=i_\mu^\ast\omega_B$, $h_\mu\cdot \pi_\mu=H\cdot i_\mu$,
$F(\mathbf{J}^{-1}(\mu))\subset \mathbf{J}^{-1}(\mu) $, $f_\mu\cdot
\pi_\mu=\pi_\mu \cdot F\cdot i_\mu$, $\mathcal{C}\cap
\mathbf{J}^{-1}(\mu)\neq \emptyset $, and
$\mathcal{C}_\mu=\pi_\mu(\mathcal{C}\cap \mathbf{J}^{-1}(\mu))$, is
a CMH system, which is simply written as $R_P$-reduced
CMH system. Here $((T^\ast Q)_\mu,\omega^B_\mu)$ is the
$R_P$-reduced magnetic symplectic space, the function $h_\mu:(T^\ast
Q)_\mu\rightarrow \mathbb{R}$ is the reduced Hamiltonian, the
fiber-preserving map $f_\mu:(T^\ast Q)_\mu\rightarrow (T^\ast
Q)_\mu$ is the reduced (external) force map, $\mathcal{C}_\mu$ is a
fiber submanifold of \;$(T^\ast Q)_\mu$ and is the reduced control
subset.
\end{defi}

It is worthy of noting that for the regular point reducible CMH
system $(T^\ast Q,G,\omega_B,H,F,\mathcal{C})$, the
$G$-invariant external force map $F: T^*Q \rightarrow T^*Q $ has to
satisfy the conditions $F(\mathbf{J}^{-1}(\mu))\subset
\mathbf{J}^{-1}(\mu), $ and $f_\mu\cdot \pi_\mu=\pi_\mu \cdot F\cdot
i_\mu, $ such that we can define the reduced external force map
$f_\mu:(T^\ast Q)_\mu\rightarrow (T^\ast Q)_\mu. $ The condition
$\mathcal{C} \cap \mathbf{J}^{-1}(\mu)\neq \emptyset $ in above
definition makes that the $G$-invariant control subset
$\mathcal{C}\cap \mathbf{J}^{-1}(\mu)$ can be reduced and the
reduced control subset is $\mathcal{C}_\mu=\pi_\mu(\mathcal{C}\cap
\mathbf{J}^{-1}(\mu))$. If the control subset cannot be reduced, we
cannot get the $R_P$-reduced CMH system.\\

From the above discussion, we know that the magnetic term is an
important global notion in the phase space of a CMH system.
it is different from the regular point reduction of a RCH system
defined on a cotangent bundle with the canonical structure,
the regular point reduction of a CMH system with symmetry reveals
the deeper relationship of the intrinsic geometrical structure
of phase space of the CMH system.

\section{CMH System with Symmetry of the
Heisenberg Group}

In this section, we first define a CMH system with symmetry of the
Heisenberg group $\mathcal{H} $,
by using the magnetic term on the cotangent bundle
of the Heisenberg group. Then we state that the CMH
system with symmetry of the Heisenberg group
is a regular point reducible CMH system, and
give its regular point reduced CMH
system on the generalization of coadjoint orbit of the
Heisenberg group $\mathcal{H} $ by calculation in detail.
Moreover, we discuss the MR-CH-equivalence for the
regular point reducible CMH system with symmetry of
the Heisenberg group $\mathcal{H} $, and prove the regular point reduction
theorem for such system.\\

A CMH system
with symmetry of the Heisenberg group $\mathcal{H}$ is a
6-tuple $(T^\ast Q,\mathcal{H},\omega_Q,H,F,\mathcal{C})$,
where the configuration space $Q=
\mathcal{H}\times V, \; \mathcal{H}= \mathbb{R}^2\oplus \mathbb{R},
$ and $V$ is a $k$-dimensional vector space, and the cotangent
bundle $T^*Q$ with magnetic symplectic form $\omega_Q= \Omega_0-
\pi^*_Q \bar{B},$ where $\Omega_0$ is the usual canonical symplectic
form on $T^*Q$, and $\bar{B}= \pi_1^*B$ is the closed two-form on
$Q$, $B$ is a closed two-form on $\mathcal{H}$ and $\pi_1:
Q=\mathcal{H}\times V \rightarrow \mathcal{H}$, $\pi_1^*: T^*
\mathcal{H}\rightarrow T^*Q$ and $\pi_Q: T^* Q\rightarrow Q$,
$\pi_Q^*: T^* Q \rightarrow T^*T^* Q$, and the
Hamiltonian $H:T^\ast Q \to \mathbb{R}$, the fiber-preserving map
$F: T^\ast Q \to T^\ast Q$ and the fiber submanifold $\mathcal{C}$
of $T^\ast Q$ are all left cotangent lifted $\mathcal{H}$-action
$\Phi^{T^*}$ invariant. At first, we define the left
$\mathcal{H}$-action $\Phi$ on $Q$ as follows
\begin{equation}
\Phi: \mathcal{H}\times Q \rightarrow Q, \;\;\;
\Phi((u,\alpha),((v,\beta),\theta)):=((u,\alpha)(v,\beta),\theta),
\end{equation} for any $(u,\alpha),(v,\beta) \in \mathcal{H},\;
\theta \in V$, that is, the $\mathcal{H}$-action on $Q$ is the left
translation on the first factor $\mathcal{H}$, and $\mathcal{H}$
acts trivially on the second factor $V$. Because locally,
$T^\ast Q \cong T^\ast \mathcal{H} \times T^* V$, and $T^\ast V\cong V\times V^\ast$,
by using the local left trivialization of $T^* \mathcal{H}$, we have that
$T^\ast Q\cong \mathcal{H}\times \mathfrak{\eta}^\ast \times V
\times V^\ast$ (locally). For the left $\mathcal{H}$-action $\Phi:
\mathcal{H}\times Q \rightarrow Q $, the cotangent lift of the action
to its cotangent bundle $T^\ast Q$ is given by
\begin{equation}
\Phi^{T^*}: \mathcal{H}\times T^*Q \rightarrow T^*Q, \;
\Phi^{T^*}((u,\alpha),((v,\beta),(\mu,\nu),\theta,\lambda)):
=((u,\alpha)(v,\beta),(\mu,\nu),\theta,\lambda),
\end{equation} for
any $(u,\alpha),(v,\beta) \in \mathcal{H},\; (\mu,\nu) \in
\mathfrak{\eta}^\ast, \; \theta \in V, \; \lambda \in V^\ast .$
If the $\mathcal{H}$-action $\Phi$ is free and proper, then the
$\Phi^{T^*}$-action is also free and proper, and the orbit space $(T^\ast
Q)/ \mathcal{H}$ is a smooth manifold and $\pi: T^*Q \rightarrow
(T^*Q )/ \mathcal{H} $ is a smooth submersion. Since $\mathcal{H}$
acts trivially on $\mathfrak{\eta}^\ast$, $V$ and $V^\ast$, it
follows that $(T^\ast Q)/ \mathcal{H}$ is diffeomorphic to
$\mathfrak{\eta}^\ast \times V \times V^\ast$.\\

From $\S 2$ we have known that for $(\mu,\nu) \in
\mathfrak{\eta}^\ast$, the coadjoint orbit $\mathcal{O}_{(\mu,\nu)}
\subset \mathfrak{\eta}^\ast$ has the magnetic orbit symplectic
forms $\omega^\pm_{\mathcal{O}_{(\mu,\nu)}}$ given by (2.10). Let
$\omega_V$ be the canonical symplectic form on $T^\ast V \cong V
\times V^\ast$ given by $$\omega_V((\theta_1, \lambda_1),(\theta_2,
\lambda_2))=<\lambda_2,\theta_1> -<\lambda_1,\theta_2> ,$$ where
$(\theta_i, \lambda_i)\in V\times V^\ast, \; i=1,2$, $<\cdot,\cdot>$
is the natural pairing between $V^\ast$ and $V$. Thus, we can induce
the magnetic symplectic form
$\tilde{\omega}^{-}_{\tilde{\mathcal{O}}_{(\mu,\nu)}}=
\pi_{\mathcal{O}_{(\mu,\nu)}}^\ast
\omega^{-}_{\mathcal{O}_{(\mu,\nu)}}+ \pi_V^\ast \omega_V $ on
the smooth manifold
$\tilde{\mathcal{O}}_{(\mu,\nu)}=\mathcal{O}_{(\mu,\nu)} \times V
\times V^\ast ,$ where the maps $\pi_{\mathcal{O}_{(\mu,\nu)}}:
\mathcal{O}_{(\mu,\nu)} \times V \times V^\ast \to
\mathcal{O}_{(\mu,\nu)} $ and $\pi_V: \mathcal{O}_{(\mu,\nu)} \times
V \times V^\ast
\to V\times V^\ast $ are canonical projections. \\

On the other hand, from $T^\ast Q = T^\ast \mathcal{H} \times T^\ast
V$ we know that there is a canonical symplectic form $\Omega_0=
\pi^\ast_1 \omega_0 +\pi^\ast_2 \omega_V $ on $T^\ast Q$, where
$\omega_0$ is the canonical symplectic form on $T^\ast \mathcal{H}$
and the maps $\pi_1: Q= \mathcal{H}\times V \to \mathcal{H}$ and
$\pi_2: Q= \mathcal{H}\times V \to V$ are canonical projections.
Then the magnetic symplectic form on $T^* Q$ is given by $\omega_Q=
\pi^\ast_1 \omega_0 +\pi^\ast_2 \omega_V- \pi^*_Q \cdot\pi_1^*B .$
Assume that the cotangent lift of left $\mathcal{H}$-action $\Phi^{T^*}:
\mathcal{H} \times T^\ast Q \to T^\ast Q$ is symplectic with respect
to the magnetic symplectic form $\omega_Q$, and admits an associated
$\operatorname{Ad}^\ast$-equivariant momentum map $\mathbf{J}_Q:
T^\ast Q \to \mathfrak{\eta}^\ast$ such that $\mathbf{J}_Q\cdot
\pi^\ast_1=\mathbf{J}_{\mathcal{H}} ,$ where
$\mathbf{J}_{\mathcal{H}}:T^\ast \mathcal{H} \rightarrow
\mathfrak{\eta}^\ast$ is a momentum map of left $\mathcal{H}$-action
on $T^\ast \mathcal{H}$, and $\pi^\ast_1: T^\ast \mathcal{H} \to
T^\ast Q$. If $(\mu,\nu)\in\mathfrak{\eta}^\ast$ is a regular value
of $\mathbf{J}_Q$, then $(\mu,\nu)\in\mathfrak{\eta}^\ast$ is also a
regular value of $\mathbf{J}_{\mathcal{H}}$ and
$\mathbf{J}_Q^{-1}((\mu,\nu))\cong
\mathbf{J}_{\mathcal{H}}^{-1}((\mu,\nu))\times V \times V^\ast .$
Denote by $G_{(\mu,\nu)}=\{(u,\alpha)\in
\mathcal{H}|\operatorname{Ad}_{(u,\alpha)}^\ast (\mu,\nu)=(\mu,\nu)
\}$ the isotropy subgroup of coadjoint $\mathcal{H}$-action at the
point $(\mu,\nu)\in \mathfrak{\eta}^\ast$. It follows that
$G_{(\mu,\nu)}$ acts also freely and properly on
$\mathbf{J}_Q^{-1}((\mu,\nu))$, the regular point reduced space
$$(T^\ast
Q)_{(\mu,\nu)}=\mathbf{J}_Q^{-1}((\mu,\nu))/G_{(\mu,\nu)}\cong
(T^\ast \mathcal{H})_{(\mu,\nu)} \times V \times V^\ast $$ of
$(T^\ast Q,\omega_Q)$ at $(\mu,\nu)$, is a symplectic manifold with
the reduced magnetic symplectic form $\omega_{(\mu,\nu)}$ uniquely
characterized by the relation
\begin{equation}
\pi_{(\mu,\nu)}^\ast \omega_{(\mu,\nu)}=i_{(\mu,\nu)}^\ast
\omega_Q=i_{(\mu,\nu)}^\ast \pi^\ast_1 \omega_0+ i_{(\mu,\nu)}^\ast
\pi^\ast_2 \omega_V-i_{(\mu,\nu)}^\ast \pi^*_Q \cdot\pi_1^*B ,
\end{equation} where the map
$i_{(\mu,\nu)}:\mathbf{J}_Q^{-1}((\mu,\nu))\rightarrow T^\ast Q $ is
the inclusion and
$\pi_{(\mu,\nu)}:\mathbf{J}_Q^{-1}((\mu,\nu))\rightarrow (T^\ast
Q)_{(\mu,\nu)}$ is the projection. From Theorem 3.1 we have seen
that $((T^\ast \mathcal{H})_{(\mu,\nu)},\omega_{(\mu,\nu)})$ is
symplectically diffeomorphic to
$(\mathcal{O}_{(\mu,\nu)},\omega_{\mathcal{O}_{(\mu,\nu)}}^{-})$,
hence we obtain that $((T^\ast Q)_{(\mu,\nu)},\omega_{(\mu,\nu)})$
is symplectically diffeomorphic to
$\tilde{\mathcal{O}}_{(\mu,\nu)}=(\mathcal{O}_{(\mu,\nu)} \times
V\times V^\ast,\;
\tilde{\omega}_{\tilde{\mathcal{O}}_{(\mu,\nu)}}^{-}). $

\begin{rema}
If $(T^\ast Q, \omega_Q)$ is a connected magnetic symplectic
manifold, and $\mathbf{J}_Q: T^\ast Q\rightarrow
\mathfrak{\eta}^\ast$ is a non-equivariant momentum map with a
non-equivariance group one-cocycle $\sigma: \mathcal{H}\rightarrow
\mathfrak{\eta}^\ast$, which is defined by
$\sigma((u,\alpha)):=\mathbf{J}_Q((u,\alpha)\cdot
z)-\operatorname{Ad}^\ast_{(u,\alpha)^{-1}}\mathbf{J}_Q(z)$, where
$(u,\alpha)\in \mathcal{H}$ and $z\in T^\ast Q. $ Then we know that
$\sigma$ produces a new affine action $\Theta: \mathcal{H}\times
\mathfrak{\eta}^\ast \rightarrow \mathfrak{\eta}^\ast $ defined by
\begin{equation}
\Theta((u,\alpha),(\mu,\nu)):=\operatorname{Ad}^\ast_{(u,\alpha)^{-1}}(\mu,\nu)
+ \sigma((u,\alpha)),
\end{equation} where $(u,\alpha)\in
\mathcal{H}, \; (\mu,\nu) \in \mathfrak{\eta}^\ast$, with respect to
which the given momentum map $\mathbf{J}_Q$ is equivariant. Since
$\mathcal{H}$ acts freely and properly on $T^\ast Q$, and
$\tilde{G}_{(\mu,\nu)}$ denotes the isotropy subgroup of $(\mu,\nu)
\in \mathfrak{\eta}^\ast$ relative to this affine action $\Theta$
and $(\mu,\nu)$ is a regular value of $\mathbf{J}_Q$. Then the
regular point reduced space $(T^\ast
Q)_{(\mu,\nu)}=\mathbf{J}_Q^{-1}((\mu,\nu))/\tilde{G}_{(\mu,\nu)}$
is also a symplectic manifold with the reduced magnetic symplectic form
$\omega_{(\mu,\nu)}$ uniquely characterized by (5.3), and this space
is symplectically diffeomorphic to
$\tilde{\mathcal{O}}_{(\mu,\nu)}=(\mathcal{O}_{(\mu,\nu)} \times
V\times V^\ast,\;
\tilde{\omega}_{\tilde{\mathcal{O}}_{(\mu,\nu)}}^{-}), $ where
$\mathcal{O}_{(\mu,\nu)}$ is the coadjoint orbit at $(\mu,\nu) \in
\mathfrak{\eta}^\ast $ with the reduced magnetic symplectic form
$\omega_{(\mu,\nu)}$.
\end{rema}

Now assume that Hamiltonian
$H((u,\alpha),(\rho,\tau),\theta,\lambda): T^\ast Q \cong
\mathcal{H}\times \mathfrak{\eta}^\ast \times V \times V^\ast \to
\mathbb{R}$ is left cotangent lifted $\mathcal{H}$-action
$\Phi^{T^*}$ invariant, for the regular value $(\mu,\nu)\in
\mathfrak{\eta}^\ast$ of $\mathbf{J}_Q$, we have the associated
reduced Hamiltonian $h_{(\mu,\nu)}((\rho,\tau),\theta,\lambda):
(T^\ast Q)_{(\mu,\nu)} \cong\mathcal{O}_{(\mu,\nu)} \times V \times
V^\ast \to \mathbb{R}$, defined by $h_{(\mu,\nu)}\cdot
\pi_{(\mu,\nu)}=H\cdot i_{(\mu,\nu)}$, and the reduced Hamiltonian
vector field $X_{h_{(\mu,\nu)}}$ is given by the reduced magnetic Hamilton's
equation $\mathbf{i}_{X_{h_{(\mu,\nu)}}}\omega_{(\mu,\nu)}=
\mathbf{d}h_{(\mu,\nu)}. $ Moreover, assume that the
fiber-preserving map $F: T^\ast Q \to T^\ast Q$ and the fiber
submanifold $\mathcal{C}$ of $T^\ast Q$ are all left cotangent
lifted $\mathcal{H}$-action $\Phi^{T^*}$ invariant, and satisfy the
conditions in Definition 4.6, then the 6-tuple $(T^\ast
Q,\mathcal{H},\omega_Q,H,F,\mathcal{C})$ is a regular point
reducible CMH system. Thus, there exists a point
$(\mu,\nu)\in\mathfrak{\eta}^\ast$, the regular value of the
momentum map $\mathbf{J}_Q: T^\ast Q \rightarrow
\mathfrak{\eta}^\ast$, and for the given $\mathcal{H}$-invariant
external force $F: T^\ast Q \to T^\ast Q $ and $\mathcal{H}$-invariant
feedback control $u: T^*Q \rightarrow \mathcal{C}$, where
$F(\mathbf{J}_Q^{-1}(\mu,\nu))\subset \mathbf{J}_Q^{-1}(\mu,\nu) $,
and $\mathcal{C}\cap \mathbf{J}_Q^{-1}(\mu,\nu) \neq \emptyset, $
the regular point reduced CMH system is the 5-tuple
$(\mathcal{O}_{(\mu,\nu)} \times V\times
V^\ast,\tilde{\omega}_{\mathcal{O}_{(\mu,\nu)} \times V \times
V^\ast }^{-},h_{(\mu,\nu)},f_{(\mu,\nu)},u_{(\mu,\nu)})$, where
$\mathcal{O}_{(\mu,\nu)} \subset \mathfrak{\eta}^\ast$ is the
coadjoint orbit, $\tilde{\omega}_{\mathcal{O}_{(\mu,\nu)} \times V
\times V^\ast }^{-}$ is the magnetic orbit symplectic form on
$\mathcal{O}_{(\mu,\nu)} \times V\times V^\ast $,
$h_{(\mu,\nu)}\cdot \pi_{(\mu,\nu)}=H\cdot i_{(\mu,\nu)}$,
$f_{(\mu,\nu)}\cdot \pi_{(\mu,\nu)}=\pi_{(\mu,\nu)}\cdot F\cdot
i_{(\mu,\nu)}$, and $u_{(\mu,\nu)} \in
\mathcal{C}_{(\mu,\nu)}=\pi_{(\mu,\nu)}(\mathcal{C}\cap
\mathbf{J}_Q^{-1}(\mu,\nu))\subset \mathcal{O}_{(\mu,\nu)} \times V
\times V^\ast$, $u_{(\mu,\nu)}\cdot
\pi_{(\mu,\nu)}=\pi_{(\mu,\nu)}\cdot u\cdot i_{(\mu,\nu)}$.
Moreover, assume that the dynamical vector field of the
regular point reduced CMH system can be expressed by
\begin{equation}X_{(\mathcal{O}_{(\mu,\nu)} \times V\times
V^\ast,\tilde{\omega}_{\mathcal{O}_{(\mu,\nu)} \times V \times
V^\ast }^{-},h_{(\mu,\nu)},f_{(\mu,\nu)},u_{(\mu,\nu)})}
=X_{h_{(\mu,\nu)}}+
\textnormal{vlift}(f_{(\mu,\nu)})+\textnormal{vlift}(u_{(\mu,\nu)}),
\; \label{3.10}
\end{equation}
where $X_{h_{(\mu,\nu)}} \in T(\mathcal{O}_{(\mu,\nu)} \times
V\times V^*) $ is the magnetic Hamiltonian vector field of the
reduced Hamiltonian $h_{(\mu,\nu)}: \mathcal{O}_{(\mu,\nu)} \times V
\times V^\ast \to \mathbb{R}$, and
$\textnormal{vlift}(f_{(\mu,\nu)})=
\textnormal{vlift}(f_{(\mu,\nu)})X_{h_{(\mu,\nu)}} \in
T(\mathcal{O}_{(\mu,\nu)} \times V\times V^*)$,
$\textnormal{vlift}(u_{(\mu,\nu)})=
\textnormal{vlift}(u_{(\mu,\nu)})X_{h_{(\mu,\nu)}} \in
T(\mathcal{O}_{(\mu,\nu)} \times V\times V^*)$, and satisfy the
condition
\begin{equation}X_{(\mathcal{O}_{(\mu,\nu)} \times V\times
V^\ast,\tilde{\omega}_{\mathcal{O}_{(\mu,\nu)} \times V \times
V^\ast }^{-},h_{(\mu,\nu)},f_{(\mu,\nu)},u_{(\mu,\nu)})}\cdot
\pi_{(\mu,\nu)}=T\pi_{(\mu,\nu)}\cdot X_{(T^\ast
Q,\mathcal{H},\omega_Q,H,F,u)}\cdot i_{(\mu,\nu)}. \label{3.11}
\end{equation}
Note that $\mbox{vlift}(u_{(\mu,\nu)})X_{h_{(\mu,\nu)}}$ is the
vertical lift of vector field $X_{h_{(\mu,\nu)}}$ under the action
of $u_{(\mu,\nu)}$ along fiber, that is,
\begin{align*}
\mbox{vlift}(u_{(\mu,\nu)})X_{h_{(\mu,\nu)}}((\rho,\tau),\theta,\lambda)
& = \mbox{vlift}((Tu_{(\mu,\nu)}
X_{h_{(\mu,\nu)}})(u_{(\mu,\nu)}((\rho,\tau),\theta,\lambda)),
((\rho,\tau),\theta,\lambda))\\ & = (Tu_{(\mu,\nu)}
X_{h_{(\mu,\nu)}})^v_{\tilde{\sigma}} ((\rho,\tau),\theta,\lambda),
\end{align*}
where $(\mu,\nu), \; (\rho,\tau) \in \mathfrak{\eta}^\ast, \; \theta
\in V, \; \lambda \in V^\ast ,$ and $\tilde{\sigma}$ is a geodesic
in $\mathcal{O}_{(\mu,\nu)} \times V\times V^*$ connecting
$u_{(\mu,\nu)}((\rho,\tau),\theta,\lambda)$ and
$((\rho,\tau),\theta,\lambda)$, and $(Tu_{(\mu,\nu)}
X_{h_{(\mu,\nu)}})^v_{\tilde{\sigma}} ((\rho,\tau),\theta,\lambda) $
is the parallel displacement of vertical vector $(Tu_{(\mu,\nu)}
X_{h_{(\mu,\nu)}})^v ((\rho,\tau),\theta,\lambda) $ along the
geodesic $\tilde{\sigma}$ from
$u_{(\mu,\nu)}((\rho,\tau),\theta,\lambda)$ to
$((\rho,\tau),\theta,\lambda)$, and
$\mbox{vlift}(f_{(\mu,\nu)})X_{h_{(\mu,\nu)}}$ is defined in the
similar manner, see Marsden et al.\cite{mawazh10} and Theorem
4.2 in $\S 4$. In consequence, from Definition 4.6 and the above
discussion, we can get the following theorem.

\begin{theo}
The 6-tuple $(T^\ast Q,\mathcal{H},\omega_Q,H,F,\mathcal{C})$ is
a regular point reducible CMH system with symmetry of the
Heisenberg group $\mathcal{H}$, where $Q=\mathcal{H}\times V$, and
$\mathcal{H}= \mathbb{R}^2\oplus \mathbb{R}$ is the Heisenberg group
with Lie algebra $\mathfrak{\eta}\cong \mathbb{R}^2\oplus
\mathbb{R}$ and its dual $\mathfrak{\eta}^*\cong \mathbb{R}^2\oplus
\mathbb{R}$, and $V$ is a $k$-dimensional vector space, and the
Hamiltonian $H:T^\ast Q \to \mathbb{R}$, the fiber-preserving map
$F: T^\ast Q \to T^\ast Q$ and the fiber submanifold $\mathcal{C}$
of $T^\ast Q$ are all left cotangent lifted $\mathcal{H}$-action
$\Phi^{T^*}$ invariant. For a given point
$(\mu,\nu)\in\mathfrak{\eta}^\ast$, the regular value of the
momentum map $\mathbf{J}_Q: T^\ast Q \rightarrow
\mathfrak{\eta}^\ast$, and the given $\mathcal{H}$-invariant
external force $F: T^\ast Q \to T^\ast Q $ and $\mathcal{H}$-invariant
feedback control $u: T^*Q \rightarrow \mathcal{C}$, where
$F(\mathbf{J}_Q^{-1}(\mu,\nu))\subset \mathbf{J}_Q^{-1}(\mu,\nu) $,
and $\mathcal{C} \cap \mathbf{J}_Q^{-1}(\mu,\nu)\neq \emptyset $,
the regular point reduced CMH system is the 5-tuple
$(\mathcal{O}_{(\mu,\nu)} \times V\times
V^\ast,\tilde{\omega}_{\mathcal{O}_{(\mu,\nu)} \times V \times
V^\ast }^{-},h_{(\mu,\nu)},f_{(\mu,\nu)},u_{(\mu,\nu)})$, where
$\mathcal{O}_{(\mu,\nu)} \subset \mathfrak{\eta}^\ast$ is the
coadjoint orbit, $\tilde{\omega}_{\mathcal{O}_{(\mu,\nu)} \times V
\times V^\ast }^{-}$ is the magnetic orbit symplectic form on
$\mathcal{O}_{(\mu,\nu)} \times V\times V^\ast $,
$h_{(\mu,\nu)}\cdot \pi_{(\mu,\nu)}=H\cdot i_{(\mu,\nu)}$,
$f_{(\mu,\nu)}\cdot \pi_{(\mu,\nu)}=\pi_{(\mu,\nu)}\cdot F\cdot
i_{(\mu,\nu)}$, and $u_{(\mu,\nu)} \in
\mathcal{C}_{(\mu,\nu)}=\pi_{(\mu,\nu)}(\mathcal{C}\cap
\mathbf{J}_Q^{-1}(\mu,\nu)) $, $u_{(\mu,\nu)}\cdot
\pi_{(\mu,\nu)}=\pi_{(\mu,\nu)}\cdot u\cdot i_{(\mu,\nu)}$.
\end{theo}

Moreover, for a given regular point reducible CMH system
with symmetry of the Heisenberg group
$(T^\ast Q,\mathcal{H},\omega_Q,H,F,\mathcal{C})$,
assume that the dynamical vector field of the regular point reduced CMH
system $(\mathcal{O}_{(\mu,\nu)} \times V\times
V^\ast,\tilde{\omega}_{\mathcal{O}_{(\mu,\nu)} \times V \times
V^\ast }^{-},h_{(\mu,\nu)},f_{(\mu,\nu)},u_{(\mu,\nu)})$
is given by $(5.5)$, and satisfies the equation $(5.6)$.
Thus, in order to emphasize explicitly
the impact of external force and control for the reduced CMH system,
we can describe the feedback control law to modify the
structure of the regular point reducible CMH system with
symmetry of the Heisenberg group. At first, from Definition 4.5 we can give
the magnetic reducible controlled Hamiltonian matching conditions and
MR-CH-equivalence for the CMH system with symmetry of the
Heisenberg group as follows.

\begin{defi}(MR-CH-equivalence)
Suppose that we have two regular point reducible CMH
systems with symmetries of the Heisenberg group $(T^\ast Q_i,
\mathcal{H},\omega_{Q_i},H_i, F_i, \mathcal{C}_i),\; i=1,2$, we say
them to be MR-CH-equivalent, or simply, $(T^\ast Q_1,
\mathcal{H},\omega_{Q_1},H_1,F_1,\mathcal{C}_1)\stackrel{MRCH}{\sim}
(T^\ast Q_2,\mathcal{H},\omega_{Q_2},H_2,F_2,\mathcal{C}_2)$, if
there exists a diffeomorphism $\varphi:Q_1\rightarrow Q_2$ such that
the following magnetic reducible controlled Hamiltonian matching
conditions hold:

\noindent {\bf MR-CH-1:} The cotangent lifted map
$\varphi^\ast:T^\ast Q_2\rightarrow T^\ast Q_1$ is symplectic with
respect to their magnetic symplectic forms.

\noindent {\bf MR-CH-2:} For $(\mu_i,\nu_i)\in
\mathfrak{\eta}^\ast_i , \; i=1,2, $ the regular reducible points of
the CMH systems with symmetries of the Heisenberg group $(T^\ast
Q_i, \mathcal{H},\omega_{Q_i}, H_i, F_i, \mathcal{C}_i)$, $i=1,2$,
the map
$\varphi_{(\mu,\nu)}^\ast=i_{(\mu_1,\nu_1)}^{-1}\cdot\varphi^\ast\cdot
i_{(\mu_2,\nu_2)}: \mathbf{J}_{Q_2}^{-1}(\mu_2,\nu_2)\rightarrow
 \mathbf{J}_{Q_1}^{-1}(\mu_1,\nu_1)$ is
$(G_{2(\mu_2,\nu_2)},G_{1(\mu_1,\nu_1)})$-equivariant and
$\mathcal{C}_1\cap
\mathbf{J}_{Q_1}^{-1}(\mu_1,\nu_1)=\varphi_{(\mu,\nu)}^\ast
(\mathcal{C}_2\cap \mathbf{J}_{Q_2}^{-1}(\mu_2,\nu_2))$, where
$(\mu,\nu)=((\mu_1, \mu_2), (\nu_1, \nu_2))$, and denote by
$i_{(\mu_1,\nu_1)}^{-1}(S)$ the pre-image of a subset $S\subset
T^\ast Q_1$ for the map
$i_{(\mu_1,\nu_1)}:\mathbf{J}_{Q_1}^{-1}(\mu_1,\nu_1)\rightarrow
T^\ast Q_1$.

\noindent {\bf MR-CH-3:} $Im[X_{H_1}+ \textnormal{vlift}(F_1)-
T\varphi^\ast X_{H_2}-\textnormal{vlift}(\varphi^\ast
F_2\varphi_\ast)]\subset\textnormal{vlift}(\mathcal{C}_1)$, where
$Im$ means the pointwise image of the map in brackets.
\end{defi}

Next, if restricting on the set of CMH systems with symmetries of the Heisenberg group,
by using the method given in Marsden et al.\cite{mawazh10}, we
can prove the following regular point reduction theorem for the
CMH systems with symmetries of the Heisenberg group, which
explains the relationship between MR-CH-equivalence for the regular
point reducible CMH systems with symmetries of the Heisenberg
group and M-CH-equivalence for the associated regular point
reduced CMH systems.

\begin{theo}
Two regular point reducible CMH systems with symmetries of
the Heisenberg group $(T^\ast Q_i,\mathcal{H}, \omega_{Q_i},
H_i,F_i,\mathcal{C}_i)$, $i=1,2,$ are MR-CH-equivalent if and only
if the associated regular point reduced CMH systems
$(\mathcal{O}_{i(\mu_i,\nu_i)} \times V_i\times V_i^\ast,
\tilde{\omega}_{\mathcal{O}_{i(\mu_i,\nu_i)} \times V_i \times
V_i^\ast }^{-}, h_{i(\mu_i,\nu_i)}, f_{i(\mu_i,\nu_i)},
\mathcal{C}_{i(\mu_i,\nu_i)}) $, $i=1,2,$ are
M-CH-equivalent.
\end{theo}

{\bf Proof:} Assume that $(T^\ast Q_1,
\mathcal{H},\omega_{Q_1},H_1,F_1,\mathcal{C}_1)\stackrel{MRCH}{\sim}
(T^\ast Q_2,\mathcal{H},\omega_{Q_2},H_2,F_2,\mathcal{C}_2)$, then
from Definition 5.3 we know that there exists a diffeomorphism
$\varphi:Q_1\rightarrow Q_2$ such that $\varphi^\ast:T^\ast
Q_2\rightarrow T^\ast Q_1$ is symplectic with respect to their
magnetic symplectic forms, and for $(\mu_i,\nu_i)\in
\mathfrak{\eta}^\ast_i , \; i=1,2, $ the map
$\varphi_{(\mu,\nu)}^\ast=i_{(\mu_1,\nu_1)}^{-1}\cdot\varphi^\ast\cdot
i_{(\mu_2,\nu_2)}: \mathbf{J}_{Q_2}^{-1}(\mu_2,\nu_2)\rightarrow
 \mathbf{J}_{Q_1}^{-1}(\mu_1,\nu_1)$ is
$(G_{2(\mu_2,\nu_2)},G_{1(\mu_1,\nu_1)})$-equivariant, and
$\mathcal{C}_1\cap
\mathbf{J}_{Q_1}^{-1}(\mu_1,\nu_1)=\varphi_{(\mu,\nu)}^\ast
(\mathcal{C}_2\cap \mathbf{J}_{Q_2}^{-1}(\mu_2,\nu_2))$, and MR-CH-3
holds. From the following commutative Diagram-1:
\[
\begin{CD}
T^\ast Q_2 @<i_{(\mu_2,\nu_2)}<< \mathbf{J}_{Q_2}^{-1}(\mu_2,\nu_2) @>\pi_{(\mu_2,\nu_2)}>> (T^\ast Q_2)_{(\mu_2,\nu_2)}\\
@V\varphi^\ast VV @V\varphi^\ast_{(\mu,\nu)} VV @V\varphi^\ast_{(\mu,\nu)/\mathcal{H}}VV\\
T^\ast Q_1 @<i_{(\mu_1,\nu_1)}<< \mathbf{J}_{Q_1}^{-1}(\mu_1,\nu_1)
@>\pi_{(\mu_1,\nu_1)}>>(T^\ast Q_1)_{(\mu_1,\nu_1)}
\end{CD}
\]
$$\mbox{Diagram-1}$$
We can define a map $\varphi_{(\mu,\nu)/\mathcal{H}}^\ast:(T^\ast
Q_2)_{(\mu_2,\nu_2)}\rightarrow (T^\ast Q_1)_{(\mu_1,\nu_1)}$ such
that $\varphi_{(\mu,\nu)/\mathcal{H}}^\ast \cdot
\pi_{(\mu_2,\nu_2)}=\pi_{(\mu_1,\nu_1)}\cdot\varphi^\ast_{(\mu,\nu)}$.
Because $\varphi_{(\mu,\nu)}^\ast:
\mathbf{J}_{Q_2}^{-1}(\mu_2,\nu_2)\rightarrow
\mathbf{J}_{Q_1}^{-1}(\mu_1,\nu_1)$ is
$(G_{2(\mu_2,\nu_2)},G_{1(\mu_1,\nu_1)})$-equivariant,
$\varphi_{(\mu,\nu)/\mathcal{H}}^\ast$ is well-defined. We shall
show that $\varphi_{(\mu,\nu)/\mathcal{H}}^\ast$ is symplectic with
respect to the reduced magnetic symplectic forms, and
$\mathcal{C}_{1(\mu_1,\nu_1)}=\varphi_{(\mu,\nu)/\mathcal{H}}^\ast
(\mathcal{C}_{2(\mu_2,\nu_2)})$. In fact, since $\varphi^\ast:
T^\ast Q_2\rightarrow T^\ast Q_1$ is symplectic with respect to
their magnetic symplectic forms, the map
$(\varphi^\ast)^\ast:\Omega^2(T^\ast Q_1)\rightarrow \Omega^2(T^\ast
Q_2)$ satisfies $(\varphi^\ast)^\ast \omega_{Q_1}=\omega_{Q_2}$. By
(5.3), we have that
$i_{(\mu_i,\nu_i)}^\ast\omega_{Q_i}=\pi_{(\mu_i,\nu_i)}^\ast\omega_{i(\mu_i,\nu_i)},$
$i=1,2$, and from the following commutative Diagram-2,
\[
\begin{CD}
\Omega^2(T^\ast Q_1) @ >i_{(\mu_1,\nu_1)}^\ast>>
\Omega^2(\mathbf{J}_{Q_1}^{-1}(\mu_1,\nu_1)) @ <\pi_{(\mu_1,\nu_1)}^\ast<< \Omega^2((T^\ast Q_1)_{(\mu_1,\nu)})\\
@V(\varphi^\ast)^\ast VV @V(\varphi^\ast_{(\mu,\nu)})^\ast VV @V(\varphi^\ast_{(\mu,\nu)/\mathcal{H}})^\ast VV\\
\Omega^2(T^\ast Q_2) @>i_{(\mu_2,\nu_2)}^\ast>>
\Omega^2(\mathbf{J}_{Q_2}^{-1}(\mu_2,\nu_2))
@<\pi_{(\mu_2,\nu_2)}^\ast <<\Omega^2((T^\ast Q_2)_{(\mu_2,\nu_2)})
\end{CD}
\]
$$\mbox{Diagram-2}$$
we have that
\begin{align*}
\pi_{(\mu_2,\nu_2)}^\ast
\cdot(\varphi_{(\mu,\nu)/\mathcal{H}}^\ast)^\ast\omega_{1(\mu_1,\nu_1)}&
=(\varphi_{(\mu,\nu)/\mathcal{H}}^\ast\cdot
\pi_{(\mu_2,\nu_2)})^\ast\omega_{1(\mu_1,\nu_1)}\\
& =(\pi_{(\mu_1,\nu_1)}\cdot
\varphi_{(\mu,\nu)}^\ast)^\ast\omega_{1(\mu_1,\nu_1)}\\
& =(i_{(\mu_1,\nu_1)}^{-1}\cdot \varphi^\ast \cdot
i_{(\mu_2,\nu_2)})^\ast
\cdot\pi_{(\mu_1,\nu_1)}^\ast\omega_{1(\mu_1,\nu_1)}\\
& = i_{(\mu_2,\nu_2)}^\ast\cdot(\varphi^\ast)^\ast
\cdot(i_{(\mu_1,\nu_1)}^{-1})^\ast \cdot i_{(\mu_1,\nu_1)}^\ast
\omega_{Q_1}\\
& =i_{(\mu_2,\nu_2)}^\ast\cdot
(\varphi^\ast)^\ast\omega_{Q_1}\\
& =i_{(\mu_2,\nu_2)}^\ast\omega_{Q_2}
=\pi_{(\mu_2,\nu_2)}^\ast\omega_{2(\mu_2,\nu_2)}.
\end{align*}
Notice that $\pi_{(\mu_2,\nu_2)}^\ast$ is a surjective, thus,
$(\varphi_{(\mu,\nu)/\mathcal{H}}^\ast)^\ast\omega_{1(\mu_1,\nu_1)}=\omega_{2(\mu_2,\nu_2)}$.
Because by hypothesis $\mathcal{C}_i \cap
\mathbf{J}_{Q_i}^{-1}(\mu_i,\nu_i)\neq \emptyset$,
$\mathcal{C}_{i(\mu_i,\nu_i)}=\pi_{(\mu_i,\nu_i)}(\mathcal{C}_i\cap
\mathbf{J}_{Q_i}^{-1}(\mu_i,\nu_i)),\; i=1,2$ and $\mathcal{C}_1\cap
\mathbf{J}_{Q_1}^{-1}(\mu_1,\nu_1)=\varphi_{(\mu,\nu)}^\ast
(\mathcal{C}_2\cap \mathbf{J}_{Q_2}^{-1}(\mu_2,\nu_2))$, we have
that
\begin{align*}
\mathcal{C}_{1(\mu_1,\nu_1)} &
=\pi_{(\mu_1,\nu_1)}(\mathcal{C}_1\cap
\mathbf{J}_{Q_1}^{-1}(\mu_1,\nu_1))\\
& =\pi_{(\mu_1,\nu_1)}\cdot \varphi_{(\mu,\nu)}^\ast
(\mathcal{C}_2\cap \mathbf{J}_{Q_2}^{-1}(\mu_2,\nu_2))\\
& =\varphi_{(\mu,\nu)/\mathcal{H}}^\ast\cdot
\pi_{(\mu_2,\nu_2)}(\mathcal{C}_2\cap
\mathbf{J}_{Q_2}^{-1}(\mu_2,\nu_2))\\
&
=\varphi_{(\mu,\nu)/\mathcal{H}}^\ast(\mathcal{C}_{2(\mu_2,\nu_2)}).
\end{align*}
Next, from (4.3) and (5.5), we know that for $i=1,2$,
$$X_{(T^\ast Q_i, \mathcal{H}, \omega_{Q_i}, H_i, F_i, u_i)}
=X_{H_i}+\textnormal{vlift}(F_i)+\textnormal{vlift}(u_i),$$
\begin{align*}
& X_{(\mathcal{O}_{i(\mu_i,\nu_i)} \times V_i\times
V_i^\ast,\tilde{\omega}_{\mathcal{O}_{i(\mu_i,\nu_i)} \times V_i
\times V_i^\ast
}^{-},h_{i(\mu_i,\nu_i)},f_{i(\mu_i,\nu_i)},u_{i(\mu_i,\nu_i)})}\\
& =X_{h_{i(\mu_i,\nu_i)}}+
\textnormal{vlift}(f_{i(\mu_i,\nu_i)})+\textnormal{vlift}(u_{i(\mu_i,\nu_i)}),
\end{align*}
and from (5.6), we have that
\begin{align*}
& X_{(\mathcal{O}_{i(\mu_i,\nu_i)} \times V_i\times
V_i^\ast,\tilde{\omega}_{\mathcal{O}_{i(\mu_i,\nu_i)} \times V_i
\times V_i^\ast
}^{-},h_{i(\mu_i,\nu_i)},f_{i(\mu_i,\nu_i)},u_{i(\mu_i,\nu_i)})}\cdot
\pi_{(\mu_i,\nu_i)}\\
& = T\pi_{(\mu_i,\nu_i)}\cdot X_{(T^\ast
Q_i,\mathcal{H},\omega_{Q_i},H_i,F_i,u_i)}\cdot i_{(\mu_i,\nu_i)}.
\end{align*}
Since $H_i, F_i$ and $\mathcal{C}_i, \; i=1,2, $ are all
$\mathcal{H}$-invariant, and for $i=1,2,$
\begin{align*}
& h_{i(\mu_i,\nu_i)}\cdot \pi_{(\mu_i,\nu_i)}=H_i\cdot
i_{(\mu_i,\nu_i)},\\
&
f_{i(\mu_i,\nu_i)}\cdot\pi_{(\mu_i,\nu_i)}=\pi_{(\mu_i,\nu_i)}\cdot
F_i\cdot i_{(\mu_i,\nu_i)},\\
&
u_{i(\mu_i,\nu_i)}\cdot\pi_{(\mu_i,\nu_i)}=\pi_{(\mu_i,\nu_i)}\cdot
u_i\cdot i_{(\mu_i,\nu_i)}.
\end{align*}
From the following commutative Diagram-3,
\[
\begin{CD}
TT^\ast Q_2 @< Ti_{(\mu_2,\nu_2)}<< T\mathbf{J}_{Q_2}^{-1}(\mu_2,\nu_2) @> T\pi_{(\mu_2,\nu_2)}>> T(T^\ast Q_2)_{(\mu_2,\nu_2)}\\
@V T\varphi^\ast VV @V T\varphi^\ast_{(\mu,\nu)} VV @V T\varphi^\ast_{(\mu,\nu)/\mathcal{H}}VV\\
TT^\ast Q_1 @< Ti_{(\mu_1,\nu_1)}<<
T\mathbf{J}_{Q_1}^{-1}(\mu_1,\nu_1) @> T\pi_{(\mu_1,\nu_1)}>>
T(T^\ast Q_1)_{(\mu_1,\nu_1)}
\end{CD}
\]
$$\mbox{Diagram-3}$$
we have that
$$T\varphi_{(\mu,\nu)/\mathcal{H}}^\ast X_{h_{2(\mu_2,\nu_2)}} \cdot
\pi_{(\mu_1,\nu_1)}= T\pi_{(\mu_1,\nu_1)} \cdot T\varphi^* X_{H_2}
\cdot i_{(\mu_1,\nu_1)},$$
$$\textnormal{vlift}(\varphi_{(\mu,\nu)/\mathcal{H}}^\ast\cdot
f_{2(\mu_2,\nu_2)}\cdot \varphi_{(\mu,\nu)/\mathcal{H}\ast})\cdot
\pi_{(\mu_1,\nu_1)} =T\pi_{(\mu_1,\nu_1)}\cdot
\textnormal{vlift}(\varphi^\ast F_2\varphi_\ast)\cdot
i_{(\mu_1,\nu_1)},$$
$$\textnormal{vlift}(\varphi_{(\mu,\nu)/\mathcal{H}}^\ast\cdot
u_{2(\mu_2,\nu_2)}\cdot \varphi_{(\mu,\nu)/\mathcal{H}\ast})\cdot
\pi_{(\mu_1,\nu_1)} =T\pi_{(\mu_1,\nu_1)}\cdot
\textnormal{vlift}(\varphi^\ast u_2\varphi_\ast)\cdot
i_{(\mu_1,\nu_1)},$$ where
$\varphi_{(\mu,\nu)/\mathcal{H}\ast}=(\varphi^{-1})^\ast_{(\mu,\nu)/\mathcal{H}}:(T^\ast
Q_1)_{(\mu_1,\nu_1)}\rightarrow (T^\ast Q_2)_{(\mu_2,\nu_2)}$. From
the magnetic controlled Hamiltonian matching condition M-CH-3 we
have that
\begin{align}
& Im[(X_{h_{1(\mu_1,\nu_1)}}
+\textnormal{vlift}(f_{1(\mu_1,\nu_1)})-
T\varphi_{(\mu,\nu)/\mathcal{H}}^\ast
X_{h_{2(\mu_2,\nu_2)}} \nonumber\\
& \;\;\;\;\;\;
-\textnormal{vlift}(\varphi_{(\mu,\nu)/\mathcal{H}}^\ast \cdot
f_{2(\mu_2,\nu_2)}\cdot \varphi_{(\mu,\nu)/\mathcal{H}
\ast})]\subset \textnormal{vlift}(\mathcal{C}_{1(\mu_1,\nu_1)}).
\end{align}
So, from Definition 4.5 we get that
\begin{align*}
& (\mathcal{O}_{1(\mu_1,\nu_1)} \times V_1\times
V_1^\ast,\tilde{\omega}_{\mathcal{O}_{1(\mu_1,\nu_1)} \times V_1
\times V_1^\ast
}^{-},h_{1(\mu_1,\nu_1)},f_{1(\mu_1,\nu_1)},\mathcal{C}_{1(\mu_1,\nu_1)})\\
& \;\;\;\;\;\; \stackrel{M-CH}{\sim}(\mathcal{O}_{2(\mu_2,\nu_2)}
\times V_2\times
V_2^\ast,\tilde{\omega}_{\mathcal{O}_{2(\mu_2,\nu_2)} \times V_2
\times V_2^\ast
}^{-},h_{2(\mu_2,\nu_2)},f_{2(\mu_2,\nu_2)},\mathcal{C}_{2(\mu_2,\nu_2)}).
\end{align*}

Conversely, assume that the regular point reduced CMH systems
$(\mathcal{O}_{i(\mu_i,\nu_i)} \times V_i\times V_i^\ast, \\
\tilde{\omega}_{\mathcal{O}_{i(\mu_i,\nu_i)} \times V_i \times
V_i^\ast
}^{-},h_{i(\mu_i,\nu_i)},f_{i(\mu_i,\nu_i)},\mathcal{C}_{i(\mu_i,\nu_i)})$,
$i=1,2,$ are M-CH-equivalent. Then from Definition 4.5,
we know that there exists a diffeomorphism
$\varphi_{(\mu,\nu)/\mathcal{H}}^\ast:(T^\ast
Q_2)_{(\mu_2,\nu_2)}\rightarrow (T^\ast Q_1)_{(\mu_1,\nu_1)}$, which
is symplectic with respect to reduced magnetic symplectic forms, and
$\mathcal{C}_{1(\mu_1,\nu_1)}=\varphi_{(\mu,\nu)/\mathcal{H}}^\ast(\mathcal{C}_{2(\mu_2,\nu_2)})$,
$(\mu_i,\nu_i) \in\mathfrak{\eta}_i^\ast, \; i=1,2, $ and (5.7)
holds. We can define a map
$\varphi_{(\mu,\nu)}^\ast:\mathbf{J}^{-1}_{Q_2}(\mu_2,\nu_2)\rightarrow
\mathbf{J}^{-1}_{Q_1}(\mu_1,\nu_1)$ such that
$\pi_{(\mu_1,\nu_1)}\cdot
\varphi_{(\mu,\nu)}^\ast=\varphi_{(\mu,\nu)/\mathcal{H}}^\ast\cdot
\pi_{(\mu_2,\nu_2)}, $ and the map $\varphi^\ast: T^\ast
Q_2\rightarrow T^\ast Q_1$ such that $\varphi^\ast\cdot
i_{(\mu_2,\nu_2)}=i_{(\mu_1,\nu_1)}\cdot \varphi_{(\mu,\nu)}^\ast, $
see the commutative Diagram-1, as well as a diffeomorphism $\varphi:
Q_1\rightarrow Q_2, $ whose cotangent lift is just $\varphi^\ast:
T^\ast Q_2\rightarrow T^\ast Q_1$. From definition of
$\varphi_{(\mu,\nu)}^\ast$, we know that $\varphi_{(\mu,\nu)}^\ast$
is $(G_{2(\mu_2,\nu_2)},G_{1(\mu_1,\nu_1)})$-equivariant. In fact,
for any $z_i\in \mathbf{J}_{Q_i}^{-1}(\mu_i,\nu_i)$,
$(u_i,\alpha_i)\in G_{i(\mu_i,\nu_i)}$, $i=1,2$ such that
$z_1=\varphi_{(\mu,\nu)}^\ast(z_2)$, and
$[z_1]=\varphi^\ast_{(\mu,\nu)/\mathcal{H}}[z_2]$, then we have that
\begin{align*} \pi_{(\mu_1,\nu_1)}\cdot\varphi_{(\mu,\nu)}^\ast(\Phi_{2(u_2,\alpha_2)}(z_2))
& =\pi_{(\mu_1,\nu_1)}\cdot \varphi_{(\mu,\nu)}^\ast((u_2,\alpha_2)
\cdot z_2)\\
& =\varphi_{(\mu,\nu)/\mathcal{H}}^\ast\cdot \pi_{(\mu_2,\nu_2)}((u_2,\alpha_2) \cdot z_2)\\
& =\varphi_{(\mu,\nu)/\mathcal{H}}^\ast[z_2]=[z_1]
=\pi_{(\mu_1,\nu_1)}((u_1,\alpha_1) \cdot z_1)
\\ & =\pi_{(\mu_1,\nu_1)}(\Phi_{1(u_1,\alpha_1)}(z_1)) =\pi_{(\mu_1,\nu_1)}\cdot \Phi_{1(u_1,\alpha_1)}\cdot
\varphi_{(\mu,\nu)}^\ast(z_2).
\end{align*}
Since $\pi_{(\mu_1,\nu_1)}$ is surjective, hence,
$\varphi_{(\mu,\nu)}^\ast\cdot
\Phi_{2(u_2,\alpha_2)}=\Phi_{1(u_1,\alpha_1)}\cdot
\varphi_{(\mu,\nu)}^\ast$. Moreover, we have
\begin{align*} & \pi_{(\mu_1,\nu_1)}(\mathcal{C}_1\cap
\mathbf{J}_{Q_1}^{-1}(\mu_1,\nu_1))=\mathcal{C}_{1(\mu_1,\nu_1)}\\
&
=\varphi_{(\mu,\nu)/\mathcal{H}}^\ast(\mathcal{C}_{2(\mu_2,\nu_2)})\\
& =\varphi_{(\mu,\nu)/\mathcal{H}}^\ast\cdot
\pi_{2(\mu_2,\nu_2)}(\mathcal{C}_2\cap
\mathbf{J}_{Q_2}^{-1}(\mu_2,\nu_2))\\
& =\pi_{(\mu_1,\nu_1)}\cdot \varphi_{(\mu,\nu)}^\ast
(\mathcal{C}_2\cap \mathbf{J}_{Q_2}^{-1}(\mu_2,\nu_2)),
\end{align*} since
$\mathcal{C}_i \cap \mathbf{J}_{Q_i}^{-1}(\mu_i,\nu_i)\neq
\emptyset,\; i=1,2 $ and $\pi_{(\mu_1,\nu_1)}$ is surjective, then
we get that $\mathcal{C}_1 \cap
\mathbf{J}_{Q_1}^{-1}(\mu_1,\nu_1)=\varphi_{(\mu,\nu)}^\ast
(\mathcal{C}_2\cap \mathbf{J}_{Q_2}^{-1}(\mu_2,\nu_2))$. We shall
show that $\varphi^\ast$ is symplectic with respect to magnetic
symplectic forms. Because
$\varphi_{(\mu,\nu)/\mathcal{H}}^\ast:(T^\ast
Q_2)_{(\mu_2,\nu_2)}\rightarrow (T^\ast Q_1)_{(\mu_1,\nu_1)}$ is
symplectic with respect to reduced magnetic symplectic forms, the
map $(\varphi_{(\mu,\nu)/\mathcal{H}}^\ast)^\ast:\Omega^2((T^\ast
Q_1)_{(\mu_1,\nu_1)})\rightarrow \Omega^2((T^\ast
Q_2)_{(\mu_2,\nu_2)})$ satisfies
$(\varphi_{(\mu,\nu)/\mathcal{H}}^\ast)^\ast
\omega_{1(\mu_1,\nu_1)}=\omega_{2(\mu_2,\nu_2)}$. By (5.3), we have
that $i_{(\mu_i,\nu_i)}^\ast
\omega_{Q_i}=\pi_{(\mu_i,\nu_i)}^\ast\omega_{i(\mu_i,\nu_i)}, $
$i=1,2$, and from the commutative Diagram-2, we have that
\begin{align*}
i_{(\mu_2,\nu_2)}^\ast\omega_{Q_2} & =
\pi_{(\mu_2,\nu_2)}^\ast\omega_{2(\mu_2,\nu_2)}
=\pi_{(\mu_2,\nu_2)}^\ast\cdot(\varphi_{(\mu,\nu)/\mathcal{H}}^\ast)^\ast\omega_{1(\mu_1,\nu_1)}\\
& =(\varphi_{(\mu,\nu)/\mathcal{H}}^\ast\cdot
\pi_{(\mu_2,\nu_2)})^\ast\omega_{1(\mu_1,\nu_1)}
=(\pi_{(\mu_1,\nu_1)}\cdot \varphi_{(\mu,\nu)}^\ast)^\ast\omega_{1(\mu_1,\nu_1)}\\
& =(i_{(\mu_1,\nu_1)}^{-1}\cdot \varphi^\ast\cdot
i_{(\mu_2,\nu_2)})^\ast\cdot \pi_{(\mu_1,\nu_1)}^\ast
\omega_{1(\mu_1,\nu_1)}\\
& =i_{(\mu_2,\nu_2)}^\ast\cdot(\varphi^\ast)^\ast
\cdot(i^{-1}_{(\mu_1,\nu_1)})^\ast\cdot
i_{(\mu_1,\nu_1)}^\ast\omega_{Q_1}=i_{(\mu_2,\nu_2)}^\ast\cdot
(\varphi^\ast)^\ast\omega_{Q_1}.
\end{align*}
Notice that $i_{(\mu_2,\nu_2)}^\ast$ is injective, thus,
$\omega_{Q_2}=(\varphi^\ast)^\ast\omega_{Q_1}. $ Since the vector
field $X_{(T^\ast Q_i,\mathcal{H},\omega_{Q_i},H_i,F_i,u_i)}$ and
$X_{(\mathcal{O}_{i(\mu_i,\nu_i)} \times V_i\times
V_i^\ast,\tilde{\omega}_{\mathcal{O}_{i(\mu_i,\nu_i)} \times V_i
\times V_i^\ast
}^{-},h_{i(\mu_i,\nu_i)},f_{i(\mu_i,\nu_i)},u_{i(\mu_i,\nu_i)})}$ is
$\pi_{(\mu_i,\nu_i)}$-related, $i=1,2,$ and $H_i, F_i$ and
$\mathcal{C}_i, \; i=1,2, $ are all $\mathcal{H}$-invariant, in the
same way, from (5.7), we have that
$$Im[X_{H_1}+ \textnormal{vlift}(F_1)-
T\varphi^\ast X_{H_2}-\textnormal{vlift}(\varphi^\ast
F_2\varphi_\ast)]\subset\textnormal{vlift}(\mathcal{C}_1),$$ that
is, the magnetic reducible controlled Hamiltonian matching condition MR-CH-3
holds. Thus, from Definition 5.3 we get that
$$(T^\ast Q_1,
\mathcal{H},\omega_{Q_1},H_1,F_1,\mathcal{C}_1)\stackrel{MRCH}{\sim}
(T^\ast Q_2,\mathcal{H},\omega_{Q_2},H_2,F_2,\mathcal{C}_2). \hskip
1cm \blacksquare $$

\section{Application: The Heisenberg Particle in a Magnetic Field}

In this section, we consider the motion of a particle of mass $m$
and charge $e$ moving in the Heisenberg group $\mathcal{H}$ under
the influence of a given magnetic field $B$, where $B$ is a closed
two-form on $\mathcal{H}$. The phase space of motion of the particle
is the cotangent bundle $T^* \mathcal{H}$, which is trivialized
locally as $\mathcal{H}\times \eta^*$ with the cotangent coordinates
$(q^i, p_i), \; i=1,2,3.$ The expressions of canonical symplectic
form $\omega_0, $ the closed two-form $B$ and the magnetic
symplectic form $\omega_B$ on $T^* \mathcal{H}$ are given by
$$
\omega_0=\sum^3_{i=1} \mathbf{d}q^i \wedge \mathbf{d}p_i ,
\;\;\;\;\;\; B=\sum^3_{i,j=1}B_{ij}\mathbf{d}q^i \wedge
\mathbf{d}q^j ,\;\;\; \mathbf{d}B=0,$$
$$\omega_B= \omega_0 -\pi^*B=\sum^3_{i=1} \mathbf{d}q^i \wedge
\mathbf{d}p_i- \sum^3_{i,j=1}B_{ij}\mathbf{d}q^i \wedge
\mathbf{d}q^j.
$$
Here $q=(q^1,q^2,q^3)\in \mathcal{H}$ is the position of the
particle in $\mathcal{H}$, and $p=(p_1,p_2,p_3)\in \eta^*$ is the
momentum of the particle. Assume that there is a left-invariant
metric $<,>_{\mathcal{H}}$ on the Heisenberg group $\mathcal{H}$.
The Hamiltonian $H: T^* \mathcal{H}\rightarrow \mathbb{R}$ is given
by the kinetic energy of the particle, that is, $$H(q,p)=
\frac{1}{2m}<p,p>_{\mathcal{H}}. $$ Note that the Hamiltonian does
not dependent on the variable $q$ and hence $\frac{\partial
H}{\partial q^i}=0, \; i=1,2,3.$ From the magnetic Hamilton's
equation $\mathbf{i}_{X_H}\omega_B= \mathbf{d}H, $ we can get the
Hamiltonian vector field as follows
$$
X_H= \sum^3_{i=1} \frac{\partial H}{\partial
p_i}\frac{\partial}{\partial q^i} +\frac{2e}{c}
\sum^3_{i,j=1}B_{ij}\frac{\partial H}{\partial
p_j}\frac{\partial}{\partial p^i},
$$
where $c$ is the speed of light, and hence we obtain the
equation of motion for the Heisenberg particle.\\

Moreover, we consider the magnetic potential $A: \mathcal{H}
\rightarrow \eta^* $, which is an one-form on the Heisenberg group
$\mathcal{H}$ and $B=\mathbf{d}A. $ Then the fiber translation map
$t_A: T^* \mathcal{H}\rightarrow T^* \mathcal{H}, \;
(q,p)\rightarrow (q,p+\frac{e}{c}A)$ can pull back the canonical
symplectic form $\omega_0$ of $T^* \mathcal{H}$ to the magnetic
symplectic form $\omega_B, $ that is, $t^*_A \omega_0= \omega_0-
\pi^*\mathbf{d}A=\omega_0- \pi^*B=\omega_B ,$ where $\pi: T^*
\mathcal{H}\rightarrow \mathcal{H}$ is the natural projection. The
modified Hamiltonian $H_A: T^* \mathcal{H}\rightarrow \mathbb{R}$ is
given by $$H_A(q, p-\frac{e}{c}A)=
\frac{1}{2m}<p-\frac{e}{c}A,\;p-\frac{e}{c}A>_{\mathcal{H}}.$$ From
the canonical Hamilton's equation $\mathbf{i}_{X_{H_A}}\omega_0=
\mathbf{d}H_A, $ we can get the same Hamiltonian vector field, that
is, $X_{H_A}=X_H. $ In fact, from Marsden and Ratiu \cite{mara99} we
know why this is a general phenomenon by using the momentum shifting
lemma.\\

On the other hand, we can also consider the magnetic term from the
viewpoint of Kaluza-Klein construction. Assume that there is a
Riemannian metric $<,>_Q$ on manifold $Q=\mathcal{H}\times S^1, $
which is obtained by keeping the left-invariant metric
$<,>_{\mathcal{H}}$ on $\mathcal{H}$ and the standard metric on
$S^1$ and declaring $\mathcal{H}$ and $S^1$ orthogonal. The metric
is called the Kaluza-Klein metric on $Q$. Note that the reduced
Hamiltonian system is not the geodesic flow of the left-invariant
metric $<,>_{\mathcal{H}}$, because of the presence of the magnetic
term. However, the equation of motion of the Heisenberg particle in
the magnetic field can be obtained by Legendre transformation and
the reducing the geodesic flow of the Kaluza-Klein metric on $Q=
\mathcal{H}\times S^1$. In the following we shall state how the
magnetic term in the magnetic symplectic form $\omega_B=
\omega_0-\pi^*B$
is obtained by reduction from the Kaluza-Klein construction.\\

Assume that $Q= \mathcal{H}\times S^1$ with Lie group $G=S^1$ acting
on $Q$, which only acts on the second factor. Since the
infinitesimal generator of this action defined by $\xi \in
\mathfrak{g}\cong \mathbb{R}= \textrm{Lie}(S^1)$ has the expression
$\xi_Q(q,\theta)=(q,\theta,0,\xi), $ by using the left local
trivialization of $T^* Q$, that is, $T^* Q\cong \mathcal{H}\times
S^1 \times \mathfrak{\eta}^*\times \mathbb{R} $ (locally), the
momentum map $\mathbf{J}_Q: T^* Q\cong \mathcal{H}\times S^1 \times
\mathfrak{\eta}^* \times \mathbb{R} \rightarrow \mathfrak{g}^*
(\cong\mathbb{R})$ is given locally by $\mathbf{J}_Q(q,\theta,p,\lambda)\xi=
(p,\lambda)\cdot(0,\xi)=\lambda\xi, $ that is,
$\mathbf{J}_Q(q,\theta,p,\lambda)=\lambda. $ In this case, the
coadjoint action is trivial. For any $\mu\in \mathfrak{g}^*
(\cong\mathbb{R}), $ we have that the isotropy group $G_\mu=S^1,$ and
its Lie algebra $\mathfrak{g}_\mu=\mathbb{R}, $ and the one-form on
$Q$, $\alpha_\mu=\lambda(A_{q_1}\mathbf{d}q_1+ A_{q_2}\mathbf{d}q_2
+A_{q_3}\mathbf{d}q_3 + \mathbf{d}\theta), $ where
$\mathbf{d}\theta$ denotes the length one 1-form on $S^1$. Note that
$\alpha_\mu$ is $S^1$-invariant and its values are in
$\mathbf{J}_Q^{-1}(\mu)=\{(q,\theta,p,\lambda)\in T^* Q \; | \; q\in
\mathcal{H}, \; \theta \in S^1, \; p\in \mathfrak{\eta}^*, \;
\lambda\in \mathbb{R} \}, $ and the exterior differential of
$\alpha_\mu$ equals $\beta_\mu = \mathbf{d}\alpha_\mu=
\mu\mathbf{d}A= \mu B. $ Thus, the closed 2-form $\beta_\mu$ on the
base $Q_\mu=Q/G_\mu=(\mathcal{H}\times S^1)/S^1=\mathcal{H}, $
equals $\mu B$ and hence the magnetic term, that is, the closed
2-form $B_\mu= \pi_{Q_\mu}^* \beta_\mu, $ is also $\mu B$, since the
map $\pi_{Q_\mu}: T^* Q_\mu= \mathcal{H}\times \mathfrak{\eta}^*
\rightarrow Q_\mu=\mathcal{H} $ is the canonical projection.
Therefore, from the cotangent bundle reduction theorem---embedding
version, we know that the reduced space $((T^* Q)_\mu, \omega_\mu)$
is symplectically diffeomorphic to $(T^* \mathcal{H},
\omega_B=\omega_0-\mu B), $ which coincides with the phase space of
Hamiltonian formulation of the Heisenberg particle in a magnetic
field $B$. If we take that $\mu=e/c, $ then the magnetic term in the
magnetic symplectic form $\omega_B$ is the magnetic field $B$ up to
a factor.\\

\noindent {\bf Acknowledgments:} The author would like to thank
Professor Tudor S. Ratiu for his understanding, support and help
in the study of geometric mechanics and cooperation, and to
dedicate the article to his on the occasion of his 65th birthday.
H. Wang's research was partially supported by Nankai
University, 985 Project and the Key Laboratory of Pure Mathematics
and Combinatorics, Ministry of Education, China.\\

\end{document}